\documentclass[final]{siamltex}
\usepackage{amssymb}
\usepackage{euscript,graphicx,epsfig,amsfonts,mathtools,epstopdf}%,bbm}
\usepackage[ruled,lined,vlined,linesnumbered,algosection]{algorithm2e}

\def\cal{\mathcal}
\def\norm#1{\left\|#1\right\|}

\def\C{\mathbb{C}}

\def\R{\mathbb{R}}

\def\K{\mathbb{K}}

\def\Cn{\C^n}

\def\Cnn{\C^{n\times n}}

\def\BA{{\bf A}}  
  
  \def\CC{{\cal C}}

  \def\CK{{\cal K}}

  \def\CP{{\cal P}}
  
  \def\CR{{\cal R}}
  
  \def\CT{{\cal T}}
  \def\CU{{\cal U}}

\def\BAs{\BA{\kern-1.5pt}}

\def\CPs{\CP{\kern-0.8pt}}

\mathcode`\@="8000
{\catcode`\@=\active \gdef@{\mkern1mu}}
\def\mydate{\number\day\ {\ifcase\month \or January\or February\or
              March\or April\or May\or June\or July\or August\or
              September\or October\or November\or December\fi}
\number\year}
\def\vek#1{\mathbf{#1}}

\newcommand{\argmin}[1]{\underset{#1}{\text{{\rm argmin}}}}

\newcommand{\prn}[1]{\left(#1\right)}

\newcommand{\brac}[1]{\left[#1\right]}
\newcommand{\curl}[1]{\left\{#1\right\}}

\newcommand\restr[2]{{% we make the whole thing an ordinary symbol
  \left.\kern-\nulldelimiterspace % automatically resize the bar with \right
  #1 % the function
  \vphantom{\big|} % pretend it's a little taller at normal size
  \right|_{#2} % this is the delimiter
  }}

\def\be{\begin{equation}}
\def\ee{\end{equation}}
\def\bea{\begin{eqnarray}}
\def\eea{\end{eqnarray}}
\def\nn{\nonumber}
\def\mand{\mbox{\ \ \ and\ \ \ }}

\def\mwith{\mbox{\ \ \ with\ \ \ }}
\def\mwhere{\mbox{\ \ \ where\ \ \ }}

\def\bbmat{\begin{bmatrix}}
\def\ebmat{\end{bmatrix}}

\def\balg{\begin{algorithm}}
\def\ealg{\end{algorithm}}
\def\bthm{\begin{theorem}}
\def\ethm{\end{theorem}}
\def\blem{\begin{lemma}}
\def\elem{\end{lemma}}
\def\bprop{\begin{proposition}}
\def\eprop{\end{proposition}}
\def\bcor{\begin{corollary}}
\def\ecor{\end{corollary}}
\def\bdefin{\begin{definition}}
\def\edefin{\end{definition}}
\def\bc{\begin{cases}}
\def\ec{\end{cases}}
\def\bproof{\par\addvspace{1ex} \indent\textit{Proof}.\ \ }
\def\eproof{\hfill\cvd\linebreak\indent}

\def\bex{\begin{example}}
\def\eex{\end{example}}
\newtheorem{assumption}{Assumption}
\def\bass{\begin{assumption}}
\def\eass{\end{assumption}}

\def\cvd{~\vbox{\hrule\hbox{%
  \vrule height1.3ex\hskip0.8ex\vrule}\hrule } }

\newcommand*{\Resize}[2]{\resizebox{#1}{!}{$#2$}}%
\newcommand*{\Sf}[1]{\Resize{.03\textwidth}{\left(#1\right)}}
\def\Sfs{\Sf{\sigma}}

\title{Block Krylov subspace recycling for shifted systems with unrelated right-
hand sides\thanks{% 
This version dated \today.}}

\author{Kirk M. Soodhalter\footnotemark[2]}

\begin{document}
\maketitle
\renewcommand{\thefootnote}{\fnsymbol{footnote}}
\footnotetext[2]{Industrial Mathematics Institute, Johannes Kepler University, 
Altenbergerstra{\ss}e 69, A-4040 Linz, Austria.
({\tt kirk.soodhalter@indmath.uni-linz.ac.at})}
 
\begin{abstract}
Many Krylov subspace methods for shifted linear systems take
 advantage of the
invariance of the Krylov subspace under a shift of the matrix.
However, exploiting this fact in the non-Hermitian case introduces restrictions; e.g., initial residuals 
must be collinear and this collinearity
must be maintained at restart.  Thus 
we cannot simultaneously solve shifted systems with
unrelated right-hand sides using this strategy, 
and all shifted residuals cannot be simultaneously
minimized over a Krylov subspace such that collinearity is maintained. 
It has been shown that this renders them generally incompatible with techniques of 
subspace recycling [Soodhalter et al. APNUM '14].

This problem, however, can be overcome.  
By interpreting a family of shifted systems as one Sylvester equation,
we can take advantage of the known ``shift invariance'' of the Krylov subspace
generated by the Sylvester operator.  
Thus we can simultaneously solve all systems over one block Krylov subspace using FOM
or GMRES type methods, even when they have unrelated right-hand sides.  
Because residual collinearity
is no longer a requirement at restart, these methods are fully compatible with subspace
recycling techniques.
Furthermore, we realize the benefits of block sparse matrix operations which arise
in the context of high-performance computing applications.

In this paper,
we discuss exploiting this Sylvester equation point of view which has yielded methods for shifted
systems which are compatible with unrelated right-hand sides.
From this, we propose a recycled GMRES method for simultaneous solution of shifted systems. 
 Numerical experiments demonstrate the effectiveness of 
the methods.
\end{abstract}


%\begin{AMS}
%65F10, 65N12, 15B57, 45B05, 45A05
%\end{AMS}

\pagestyle{myheadings}
\thispagestyle{plain}
\markboth{K. M. SOODHALTER}{BLOCK KRYLOV METHODS FOR SOLVING SHIFTED SYSTEMS}

\section{Introduction}\label{section.intro}
For a given coefficient matrix $\vek A\in\C^{n\times n}$, a problem which often arises
in applied mathematics is to solve multiple linear systems in which the coefficient
matrix of each system differs from $\vek A$ by a scalar multiple of the identity, i.e.,
we must solve
\be\label{eqn.Asigxb}
(\vek A + \sigma_{i}\vek I)\vek x\Sf{\sigma_{i}} = \vek b\Sf{\sigma_{i}}\mwith i=1, 2,\ldots, L\mwhere \curl{\sigma_{i}}_{i=1}^{L}\subset\C.
\ee
Such families arise in applications such as Tikhonov-Phillips regularization, lattice quantum chromodynamics, 
rational Krylov subspaces, diffuse optical tomography, etc.  
When the coefficient matrix is large and sparse, matrix-free iterative methods such
as Krylov subspace methods are of interest.  When we
are solving multiple shifted linear systems, such as those in \eqref{eqn.Asigxb}, Krylov 
subspace methods are particularly attractive because, under certain assumptions,
the Krylov subspace $\CK_{j}(\vek A,\vek u)$, cf., \eqref{eqn.Krylov-basis}, 
is invariant under scalar shift of the coefficient matrix.  Specifically, we have that
\be\label{eqn.collinearity}
	\CK_{j}(\vek A + \sigma_{i_{1}}\vek I,\vek u) = \CK_{j}(\vek A + \sigma_{i_{2}}\vek I,\widetilde{\vek u})
\ee
as long as $\widetilde{\vek u} = \beta\vek u$ where 
$\beta \in\C\setminus\curl{0}$.  This shift invariance has led to numerous methods
for solving the systems in \eqref{eqn.Asigxb} over a single Krylov subspace;
see, e.g., 
\cite{ASV.2012, Darnell2008,Frommer2003,Frommer1998,Kilmer.deSturler.tomography.2006, Meerbergen.bai.lanczos-param-sym.2010,SBK.2014,Simoncini2003a}.  
However, this \emph{collinearity requirement} means, in general, we cannot use 
the invariance \eqref{eqn.collinearity} when the right-hand 
sides $\curl{\vek b\Sf{\sigma_{i}}}_{i=1}^{L}$ are unrelated; and 
in any case, for GMRES type methods, we 
cannot simultaneously minimize all residuals while maintaining collinearity
\cite{Frommer1998}.  Furthermore, it has been shown that the collinearity requirement causes great difficulty
when incorporating shifted system solvers into the subspace recycling framework \cite{SSX.2014}.

Therefore, we propose an alternative.  In this paper, we recall that we can 
still exploit the shift invariance of the Krylov subspace but avoid the collinearity
restriction by exploiting the invariance of a block Krylov subspace generated by a related
Sylvester operator \cite{GJR.2002,RS.2002,S.2015-website,S.1996}.
By collecting initial residuals for all systems in \eqref{eqn.Asigxb} as columns of
$\vek R_{0}^{\boldsymbol{\sigma}}\in\C^{n\times L}$ and building a block Krylov subspace, we can construct
approximations for each shifted system over one block Krylov subspace according to
a Petrov-Galerkin condition on each residual (e.g., GMRES \cite{Saad.GMRES.1986} or FOM \cite{S1981}).  By building the
block Krylov subspace from all the residuals we avoid 
the above discussed problems arising from a
lack of collinearity.  
Building upon block Krylov subspace technology allows us to use existing,
well-tested implementation strategies with minor modifications.
By avoiding the restrictive collinearity requirement, this produces methods which 
can be incorporated into the subspace recycling framework \cite{Parks.deSturler.GCRODR.2005}.
Furthermore, building methods from block Krylov
subspace techniques allows us to realize the benefits in communication
efficiency which have been observed for block Krylov methods with their 
sparse block operations; see, e.g., \cite{CK.2010,H.2010,MHDY.2009,PSS.2013}.

For clarity, it should be noted, the methods presented in this paper are \textbf{NOT}
extensions of the shifted GMRES method \cite{Frommer1998} or the shifted 
FOM \cite{Simoncini2003a} methods to block Krylov subspace to solve
problems with multiple right-hand sides.  Extensions of these methods to the case of multiple
right-hand sides can be derived and indeed do already exist; see, e.g., \cite{Darnell2008,WPS.2015}.
However, such methods still require that the columns of the block residuals for each shifted
system span the same subspace (a generalization of the collinearity requirement).  Thus they are not
applicable to problems with unrelated right-hand sides.  Here, we are treating the problem of solving
multiple shifted systems by taking advantage of its equivalence to a Sylvester equation and using
what has already been proven in that context (see, e.g., \cite{RS.2002}) to build solvers 
satisfying the requirements necessary for compatibility with the subspace recycling framework.

The rest of this paper is organized as follows.  In the next section, we review
methods for shifted linear systems based upon the invariance \eqref{eqn.collinearity} and discussed
their restrictions.  
In Section
\ref{section.bl-sh-arnoldi}, we review the Sylvester equation formulation of a family of
shifted systems and show how through an associated block Krylov subspace
invariance, we get a shifted Arnoldi relation for each individual shifted 
system in \eqref{eqn.Asigxb}, we further discuss how the block shift invariance
leads to GMRES- and FOM-like methods which allow simultaneous projection 
of all residuals according to a Petrov-Galerkin condition over the block
Krylov subspace from which we can build subspace recycling methods.  
We also discuss strategies when initial residuals are not linearly independent.
In Section \ref{section.shifted-recycling}, we propose a recycled GMRES method
for simultaneous solution of multiple shifted systems.
Further algorithmic details are also discussed.
In Section \ref{section.perf}, we discuss performance 
of these new Sylvester-based methods as well as how one
selects recycled augmentation subspaces.
Numerical results demonstrating proof of concept and effectiveness of these methods are 
presented in Section \ref{section.num-tests}.
\section{Preliminaries}\label{section.prelim}
We begin with a brief review of Krylov subspace methods as well as
techniques for solving
shifted linear system and of subspace recycling techniques.
Recall that in many Krylov subspace iterative methods for solving
the unshifted system
\begin{equation}\label{eqn.Axb}
	\vek A\vek x = \vek b
\end{equation}
with $\vek A\in\Cnn$ , we generate an orthonormal basis for 
\begin{equation}\label{eqn.Krylov-basis}
\CK_{j}(\vek A,\vek u) = \text{span}\curl{\vek u, \vek A\vek u,\ldots, \vek A^{j-1}\vek u}
\end{equation}
with the Arnoldi process, where $\vek u$ is some starting vector.  Let $\vek V_{j}\in\C^{n\times j}$ be the matrix with orthonormal columns generated by the Arnoldi process spanning $\CK_{j}(\vek A,\vek u)$.  Then we have the Arnoldi relation
\begin{equation}\label{eqn.arnoldi-relation}
	\vek A\vek V_{j} = \vek V_{j+1}\overline{\vek H}_{j}
\end{equation}
where $\overline{\vek H}_{j}\in\C^{(j+1)\times j}$ is upper Hessenberg; 
see, e.g., \cite[Section 6.3]{Saad.Iter.Meth.Sparse.2003} and 
\cite{szyld.simoncini.survey.2007}.  Let $\vek x_{0}$ be an 
initial approximation to the solution of \eqref{eqn.Axb} 
and $\vek r_{0}=\vek b-\vek A\vek x_{0}$ be the initial residual. 
At iteration $j$, we choose 
$\vek x_{j}=\vek x_{0}+\vek t_{j}$, with $\vek t_{j}\in\CK_{j}(\vek A,\vek r_{0})$.  
The vector $\vek t_{j}$ is called a \emph{correction} and the subspace from which
it is drawn (in this case a Krylov subspace) is called a \emph{search space}. 
In GMRES \cite{Saad.GMRES.1986}, $\vek t_{j}$ satisfies
\be\nn
	\vek b - \vek A(\vek x_{0} + \vek t_{j}) \perp \vek A\CK_{j}(\vek A,\vek r_{0}),
\ee
where $\vek b - \vek A(\vek x_{0} + \vek t_{j})$ is the $j$th residual which is equivalent to
\be\nn
	\vek t_{j} = \argmin{\vek t\in \CK_{j}(\vek A,\vek r_{0})}\norm{\vek b-\vek A(\vek x_{0} + \vek t)},
\ee
which is itself equivalent to solving the $(j+1)\times j$ minimization problem
\begin{equation}\label{eqn.gmres-min}
	\vek y_{j} = \argmin{\vek y\in\C^{j}}\norm{ \overline{\vek H}_{j}\vek y - \norm{\vek r_{0}}\vek e_{1}^{(j+1)}},
\end{equation}
where $\vek e_{J}^{(i)}$ denotes the $J$th Cartesian basis vector in $\C^{i}$.
We then set $\vek x_{j} = \vek x_{0} + \vek V_{j}\vek y_{j}$. Recall that in restarted 
GMRES, often called GMRES($m$), we run an $m$-step cycle of the GMRES method and 
compute an approximation $\vek x_{m}$.   We halt the process, discard $\vek V_{m}$, 
and restart with the new residual.  This process is repeated until we 
achieve convergence.  

A similar derivation leads to the related Full Orthogonalization Method (FOM) \cite{S1981}.
Here we enforce the condition that 
\be\nn
	\vek b - \vek A(\vek x_{0} + \vek t_{j}) \perp \CK_{j}(\vek A,\vek r_{0})
\ee
which is equivalent to solving the $j\times j$ linear system
\be\label{eqn.FOM-system}
	\vek H_{j}\vek y_{j} = \beta\vek e_{1}^{(j)}
\ee
where $\vek H_{j}\in\C^{j\times j}$ is simply the matrix obtained by deleting the last
row of $\overline{\vek H}_{j}$.  The iterates produced by the GMRES and FOM algorithms
are closely related; see, .e.g., \cite[Section 6.5.7]{Saad.Iter.Meth.Sparse.2003} as well as
\cite{B.1991}.
As with GMRES, a restarted version of the FOM method has been proposed called FOM($m$).

One downside to the restarting
strategy for GMRES and FOM is that one discards all the information generated in the Krylov subspace
$\CK_{j}(\vek A, \vek r_{0})$.  This leads to a characteristic delay in convergence of restarted methods 
\cite{Manteuffel.LGMRES.2005,Joubert.Convergence-Restarted-GMRES.1994}, 
as well as at times unpredictable convergence \cite{Embree.Tortoise-Restarted-GMRES.2003}.  
Krylov subspace augmentation/deflation techniques have been proposed to allow for the inclusion of vectors
to the search space in addition to the Krylov subspace 
\cite{BR-2.2007,BE.1998,deSturler.GCRO.1996,deSturler.GCROT.1999,EG.1997-tech,GGPV.FGMRES-DR.2012,Morgan.Restarted-GMRES-eig.1995,Morgan.GMRESDR.2002,Kilmer.deSturler.tomography.2006,Parks.deSturler.GCRODR.2005,SYEG.2000,WSG.2007}. 
Augmentation allows the user to select a subspace of
$\CK_{j}(\vek A,\vek r_{0})$ to include in the subsequent cycle.  Furthermore, if
one is solving a sequence of problems $\vek A_{i}\vek x_{i}=\vek b_{i}$ where 
$\Delta\vek A_{i}~=~\vek A_{i+1}~-~\vek A_{i}$ is in some sense ``not large'', one can use vectors generated
for solving system $i$ for augmentation when solving system $i+1$.  

In this paper, we focus on the augmentation
technique called subspace recycling, specifically GMRES with subspace recycling (rGMRES or GCRO-DR)
\cite{Parks.deSturler.GCRODR.2005}, and we describe this method using the framework developed by Gaul
in his thesis \cite{Gaul.2014-phd} which has overlap with related work undertaken with colleagues and
published earlier in \cite{GGL.2013}.
In principle, when one views subspace recycling methods in this framework, one sees that the 
framework can be wrapped around many iterative methods (though what results may not be easily
implemented).  Suppose we have a $k$ dimensional subspace $\CU$ of vectors which is to be part
of the search space.  Then for an appropriately chosen projector $\vek P$, we use the Krylov subspace iteration to solve the projected subproblem
\be\label{eqn.proj-prob}
	\vek P\vek A\widehat{\vek x} = \vek P\vek b.
\ee
If $\vek x_{0}$ is an initial approximation for the full problem, then  one can compute from this an initial 
approximation $\widehat{\vek x}_{0}$ for \eqref{eqn.proj-prob}.
From the iteration, we get an approximation $\widehat{\vek x}_{j} = \widehat{\vek x}_{0} + \widehat{\vek t}_{j}$, where
$\widehat{\vek t}_{j}\in\CK_{j}(\vek P\vek A,\widehat{\vek r}_{0})$.  To get an approximation for the
original problem, we compute a correction $\widehat{\vek s}_{j}\in\CU$ and set 
$\vek x_{j} = \widehat{\vek x}_{j} + \widehat{\vek s}_{j}$.  How the terms $\vek P$, $\widehat{\vek x}_{0}$, and 
$\widehat{\vek s}_{j}$ are computed is determined by $\CU$ and the iterative method being applied to the projected
problem.  

The rGMRES method \cite{Parks.deSturler.GCRODR.2005} represents the confluence of two augmentation/deflation approaches,
the harmonic Ritz vector deflation approach (GMRES-DR) of Morgan \cite{Morgan.GMRESDR.2002} and the GCRO optimal
augmentation approach of de Sturler \cite{deSturler.GCRO.1996}.  As described by Gaul \cite{Gaul.2014-phd},
this method fits into the general augmentation framework.  Let $\CC=\vek A\cdot\CU$ be the image of
$\CU$ under the action of the operator and $\vek P^{\CC}$ be the orthogonal projector onto $\CC$.  Then rGMRES proceeds
by applying a GMRES iteration to the projected subproblem
\be\label{eqn.rgmres-proj}
	\prn{\vek I - \vek P^{\CC}}\vek A\widehat{\vek x} = \prn{\vek I - \vek P^{\CC}}\vek b.
\ee
We compute the initial residual $\widehat{\vek r}_{0} = \prn{\vek I - \vek P^{\CC}}\vek r_{0}$
and  the associated initial approximation 
$\widehat{\vek x}_{0} = \vek x_{0} + \vek P_{\vek A^{\ast}\vek A}^{\CU}(\vek x_{0}-\vek x)$ which is computable without knowing
$\vek x$, where $\vek P_{\vek A^{\ast}\vek A}^{\CU}$ is the projector onto
$\CU$ orthogonal with respect to the inner product induced by $\vek A^{\ast}\vek A$; in other words,
$\vek P_{\vek A^{\ast}\vek A}^{\CU}$ projects onto the subspace $\CU$ and maps vectors in $\prn{\vek A^{\ast}\vek A\,\CU}^{\perp}$
to zero.\footnote{If $\vek U\in\C^{n\times k}$ has columns spanning $\CU$, then 
$\vek P_{\vek A^{\ast}\vek A}^{\CU}=\vek U\prn{\vek U^{T}\vek A^{\ast}\vek A\vek U}^{-1}\prn{\vek A^{\ast}\vek A\vek U}^{\ast}$.}  
The linear system \eqref{eqn.rgmres-proj} is singular and consistent.  At iteration $j$, GMRES 
produces the correction $\widehat{\vek t}_{j}$ yielding approximate solution 
$\widehat{\vek x}_{j} = \widehat{\vek x}_{0} + \widehat{\vek t}_{j}$ to \eqref{eqn.rgmres-proj}.  
To get an approximation for the full problem, we compute the correction 
$\widehat{\vek s}_{j} = \vek P_{\vek A^{\ast}\vek A}^{\CU}\widehat{\vek t}_{j}$ and set 
$\vek x_{j} = \widehat{\vek x}_{j} + \widehat{\vek s}_{j}$.  We have in this case that the residuals
for the projected problem \eqref{eqn.rgmres-proj} and the full problem \eqref{eqn.Axb} are the
same, i.e., $\vek r_{j} = \widehat{\vek r}_{j}$ and that the residual satisfies the Petrov-Galerkin
condition
\be\nn
	\vek r_{j} \perp \vek A\cdot\curl{\CU + \CK_{j}( \prn{\vek I - \vek P^{\CC}}\vek A,\widehat{\vek r}_{0} )} \subset \CC \oplus \CK_{j+1}( \prn{\vek I - \vek P^{\CC}}\vek A,\widehat{\vek r}_{0} ),
\ee
where the subspace containment arises from the modified Arnoldi relation in, e.g., \cite{Parks.deSturler.GCRODR.2005}.
Furthermore, $\widehat{\vek s}_{j}$ and $\widehat{\vek t}_{j}$ satisfy
\be\nn
	\prn{\widehat{\vek s}_{j},\widehat{\vek t}_{j}} = \argmin{\substack{\widehat{\vek s}\in\CU\\\widehat{\vek t}\in\CK_{j}( \prn{\vek I - \vek P^{\CC}}\vek A,\widehat{\vek r}_{0} )}}\norm{\vek b - \vek A(\widehat{\vek x}_{0} + \widehat{\vek s} + \widehat{\vek t})}.
\ee

Our goal here is to incorporate a strategy to simultaneously solve a family of shifted systems
into the rGMRES framework.  However, this has been shown to be difficult.
Many methods for the simultaneous solution of shifted systems 
take advantage of the shift invariance \eqref{eqn.collinearity}; see, e.g.,  
\cite{Darnell2008,Frommer2003,Frommer1998,Frommer1995,Kilmer.deSturler.tomography.2006,Kirchner2011,Simoncini2003a}.
In a non-symmetric method with restarting,
collinearity must be maintained at restart. In \cite{SSX.2014}, 
this was shown to be a
troublesome restriction when attempting to extend such techniques 
to the rGMRES setting.  In the Hermitian case, this problem is alleviated as one 
need not restart, due to the Lanczos three-term recurrence.  In this setting, 
Kilmer and de Sturler proposed a shifted recycled MINRES algorithm \cite{Kilmer.deSturler.tomography.2006}.
What we will show is that the strategies employed in \cite{Kilmer.deSturler.tomography.2006} can be used in 
the non-Hermitian case when one uses a specific type of block Krylov subspace method to solve
the family of shifted systems, such that the collinearity restriction is no longer an issue; cf., 
Section \ref{section.bl-sh-arnoldi}.

We also note that the shift-invariance \eqref{eqn.collinearity} 
no longer holds if general 
preconditioning is used.
However, specific polynomial preconditioners can be
constructed (see, e.g., 
\cite{ASV.2012,BG.2014,Jegerlehner.shifted-systems.1996, WWJ.2012}) 
for which shift invariance can be maintained.  
Since the methods proposed in this paper are built on the shift-invariance
of a Krylov subspace, they are indeed compatible with this type of 
polynomial preconditioners.
However, in this paper, we treat only the unpreconditioned problem, 
as in, e.g., \cite{Darnell2008,Frommer1998,Simoncini2003a}.

It should also be noted that methods have been proposed which do not
rely on the shift invariance property of Krylov subspace methods.  Kressner
and Tobler treated the more general situation of parameter dependent
linear systems where dependence on the parameter of the matrix and 
right-hand sides are sufficiently smooth \cite{KT.2011}. 
In \cite{S.2014}, the relationship between the shifted
coefficient matrices is exploited without using the shift invariance by solving 
one system and projecting the other residuals in a Lanczos-Galerkin procedure.

We end this section by briefly reviewing the restarted GMRES method for shifted
systems of Frommer and Gl\"assner \cite{Frommer1998}
and the restarted FOM method for shifted systems
of Simoncini \cite{Simoncini2003a}, both developed to solve
\eqref{eqn.Asigxb}.

Frommer and Gl\"{a}ssner proposed a 
restarted GMRES method to solve 
\eqref{eqn.Asigxb} in the case that the initial residuals are collinear
(which for clarity we denote in this paper as ``sGMRES \cite{Frommer1998}'').  
Within a cycle,
the residual for one system from \eqref{eqn.Asigxb} is minimized.  
We call this the {\em base residual}.
Approximations for all other systems are chosen such that their residuals
are collinear with the base residual.  This procedure reduces to solving $L-1$ 
small $(m+1)\times (m+1)$ linear systems at the end of each cycle.  Since all residuals
are then collinear at the end of the cycle, the shift invariance of the Krylov subspace
holds at the beginning of the next cycle.
It is not guaranteed for all matrices and all shifts that collinear residuals
can be computed; however, conditions are derived for when such residuals
can be constructed.
 Specifically, for a positive-real matrix 
 $\vek A$ (field of values being contained in the right half-plane), 
restarted GMRES for shifted linear systems computes solutions at every iteration for all real shifts 
$\sigma_{i} > 0$.  

Simoncini proposed an algorithm for simultaneously solving the systems in
\eqref{eqn.Asigxb} based on FOM($m$)
 (which we denote for clarity in this paper ``sFOM \cite{Simoncini2003a}'').  
 Due to the properties of the residual produced by the FOM
 algorithm, the method is conceptually simpler to describe.  For each cycle, the 
 common shift-invariant Krylov subspace is generated.  For each shifted system,
 the approximation is computed according to the Petrov-Galerkin condition which
 defines FOM.  Residuals produced by the FOM Petrov-Galerkin condition at
 step $m$
 are always collinear with the $(m+1)$st Arnoldi vector $\vek v_{m+1}$.  Therefore,
 FOM for shifted systems produces collinear residuals by default, and the Krylov
 subspace remains invariant after restart.  Thus, as long as the initial residuals
 for all shifted systems in \eqref{eqn.Asigxb} are collinear, the shifted FOM algorithm
 is applicable without modification.  However, if the right-hand sides are in general
 unrelated, we cannot use this method to simultaneously solve all linear systems in \eqref{eqn.Asigxb}.
\section{A generalization of shifted systems}\label{section.bl-sh-arnoldi}
The main tool in this paper which allows us to achieve our goals is to observe that a generalization of shifted
systems has been observed and exploited in the literature \cite{RS.2002,S.1996,S.2015-website}, and this allows one
to solve one to solve such families over block Krylov subspaces.
\subsection{Sylvester equation interpretation}
Solving the family of 
shifted systems \eqref{eqn.Asigxb} is equivalent to solving the Sylvester equation
\be\label{eqn.shifted-sylvester}
	\vek A\vek X^{\boldsymbol\sigma} - \vek X^{\boldsymbol\sigma} \vek D = \vek B^{\boldsymbol\sigma},
\ee
where $\vek D = {\rm diag}\,(\sigma_{1},\ldots,\sigma_{L})\in\C^{L\times L}$ and
$\vek X^{\boldsymbol\sigma},\vek B^{\boldsymbol\sigma}\in\C^{n\times L}$ are the block vectors 
containing the solutions and right-hand sides of each shifted system from \eqref{eqn.Asigxb}, respectively.
This observation has been previously discussed; see, e.g., \cite{RS.2002,S.1996,S.2015-website}.  Sylvester equations of
this form are in many ways a natural generalization of a single shifted system.  As shown in \cite{RS.2002},
for \emph{any} $\widetilde{\vek D}\in\C^{L\times L}$, if for $\vek F\in\C^{n\times L}$ we define the Sylvester operator as
$\CT:\vek F\mapsto \vek A\vek F + \vek F\widetilde{\vek D}$ and define powers of the operator by $\CT^{i}\vek F = \CT(\CT^{i-1}\vek F)$
with $\CT^{0}\vek F = \vek F$, then it has been shown that the block Krylov 
subspace, cf. \eqref{eqn.block-kryl-def}, generated by $\CT$ and $\vek F$ is 
equivalent to that generated by the matrix $\vek A$ and $\vek F$, i.e.,
\be\label{eqn.block-invariance}
	\K_{j}(\CT,\vek F) = \K_{j}(\vek A, \vek F).
\ee
The case in which $\widetilde{\vek D} = \vek D$ is diagonal is simply a special case of this more general result.
Let $\vek X_{0}^{\boldsymbol\sigma}$ be an initial approximate solution to \eqref{eqn.shifted-sylvester}, and
$\vek R_{0}^{\boldsymbol\sigma} = \vek B^{\boldsymbol\sigma} - \vek A\vek X_{0}^{\boldsymbol\sigma} - \vek X_{0}^{\boldsymbol\sigma}\vek D$ the associated residual.
Then one can generate the block Krylov subspace $\K_{j}(\vek A,\vek R_{0}^{\boldsymbol\sigma})$
and project \eqref{eqn.shifted-sylvester} onto $\K_{j}(\vek A,\vek R_{0}^{\boldsymbol\sigma})$
and solving a projected Sylvester equation or applying some block GMRES/FOM residual constraint.  
This is described in, e.g., \cite{GJR.2002,RS.2002,S.2015-website,S.1996}.
By exploiting the shifted system structure, we 
also can solve a FOM- or GMRES-type subproblem individually for each shifted system over
$\K_{j}(\vek A, \vek R_{0}^{\boldsymbol\sigma})$, a procedure proposed in \cite{S.1996} and described below.
In the case $\widetilde{\vek D} = \vek D$, this is equivalent to applying the FOM/GMRES methods for Sylvester 
equations to \eqref{eqn.shifted-sylvester}.

\subsection{Derivation of sbFOM and sbGMRES}
The block Krylov subspace $\K_{j}(\vek A,\vek R^{\boldsymbol\sigma}_{0})$ is a 
generalization of the definition of a Krylov subspace, i.e., 
\be\label{eqn.block-kryl-def}
	\K_{j}(\vek A,\vek R^{\boldsymbol\sigma}_{0}) = \text{span}\curl{\vek R^{\boldsymbol\sigma}_{0}, \vek A\vek R^{\boldsymbol\sigma}_{0}, \vek A^{2}\vek R^{\boldsymbol\sigma}_{0},\ldots \vek A^{j-1}\vek R^{\boldsymbol\sigma}_{0}}
\ee
where the span of a sequence of block vectors is simply the span of the columns
of all the blocks combined.
It is straightforward to show that this definition is equivalent to 
\be\label{eqn.block-krylov-sum}
	\K_{j}(\vek A,\vek R^{\boldsymbol\sigma}_{0}) = \CK_{j}(\vek A,\vek r_{0}\Sf{\sigma_{1}}) +  \CK_{j}(\vek A,\vek r_{0}\Sf{\sigma_{2}}) + \cdots +  \CK_{j}(\vek A,\vek r_{0}\Sf{\sigma_{L}}).
\ee
Except for in Section \ref{section.linear-dependence}, we assume throughout this paper
that \linebreak $\dim \K_{j}(\vek A,\vek R^{\boldsymbol\sigma}_{0}) = jL$,
which implies that the block size (i.e., the number of linearly independent right-hand sides) is
the same as the number of shifts.  Following 
the description in \cite[Section 6.12]{Saad.Iter.Meth.Sparse.2003}, we 
represent $\K_{j}(\vek A,\vek R^{\boldsymbol\sigma}_{0})$ in terms of the block Arnoldi basis 
$\curl{\vek V_{1},\vek V_{2},\ldots, \vek V_{j}}$ where 
$\vek V_{i}\in\C^{n\times L}$ has orthonormal columns and each column 
of $\vek V_{i}$ is orthogonal to all columns of $\vek V_{j}$ for all $j\neq i$.  
We obtain $\vek V_{1}$ via the reduced QR-factorization
$\vek R^{\boldsymbol\sigma}_{0}~=~\vek V_{1}\vek S_{0}$ where $\vek S_{0}\in\C^{L\times L}$ 
is upper triangular.   
Let
\be\nn
	\vek W_{j} = \brac{\begin{matrix}\vek V_{1} & \vek V_{2} \cdots \vek V_{j}\end{matrix}}\in\C^{n\times jL}.
\ee  
Let $\overline{\vek H}_{j} = (\vek H_{i\ell})\in\C^{(j+1)L\times jL}$ be the block upper Hessenberg matrix
generated by the block Arnoldi method where $\vek H_{i\ell}\in\C^{L\times L}$.  
This yields the block Arnoldi relation
\begin{equation}\label{eqn.block-arnoldi-relation}
	\vek A\vek W_{j} = \vek W_{j+1}\overline{\vek H}_{j}.
\end{equation}
A straightforward generalization of GMRES for block Krylov subspaces 
(called block GMRES) has been described \cite{Vital1990}
and a block FOM method has be similarly derived in the context of evaluating matrix exponentials \cite{WPS.2015}.

A great deal of work on the theory and implementation of 
block Krylov subspace methods has been published; see, e.g., \cite{GS.2005,GS.2009,Gutknecht2008,Simoncini.Conv-Block-GMRES,Simoncini1995,S.2014,Vital1990}.
The shift invariance properties of Krylov subspaces extend to the block
setting.  The following 
straightforward proposition directly follows from their construction.
\bprop
	The block Krylov subspace is invariant under scalar shifts
	of the coefficient matrix, i.e.,
	\be\label{eqn.block-shift-invar}
		\K_{j}(\vek A,\vek R^{\boldsymbol\sigma}_{0}) = \K_{j}(\vek A + \sigma\vek I,\vek R^{\boldsymbol\sigma}_{0})
	\ee
	 with $\sigma\in\C\setminus\curl{0}$ and satisfies the shifted 
	 block Arnoldi relation
	\be\label{eqn.blsh-arnoldi-relation}
		(\vek A + \sigma\vek I)\vek W_{j} = \vek W_{j+1}\overline{\vek H}_{j}\Sfs
	\ee
	where
	\be\nn
		\overline{\vek H}_{j}\Sfs = \overline{\vek H}_{j} + \sigma \bbmat \vek I_{jL\times jL} \\ \vek 0_{L\times jL}\ebmat.
	\ee
\eprop
The block shift invariance has
previously been exploited in \cite{Darnell2008, RS.2002, S.1996}.  
It can also be seen as a special case of \eqref{eqn.block-invariance}
with $\widetilde{\vek D} = \gamma\vek I$ for some $\gamma\in\C\setminus\curl{0}$.  
If we build the block Krylov
	subspace $\K_{j}(\vek A,\vek R_{0}^{\boldsymbol{\sigma}})$, generating the 
	associated $\vek W_{j+1}$ and $\overline{\vek H}_{j}$, then we have
	for each initial residual the two equalities
		\be\label{eqn.init-resid-relation-j}
			\vek r_{0}\Sf{\sigma_{i}} = \vek W_{j}\vek E_{1}^{(j)}\vek S_{0}\vek e_{i}^{(L)} = \vek W_{j+1}\vek E_{1}^{(j+1)}\vek S_{0}\vek e_{i}^{(L)}.
		\ee
		where $\vek R_{0}^{\boldsymbol{\sigma}} = \vek V_{1}\vek S_{0}$ is a QR-factorization as before, 
		$\vek E_{1}^{(J)}\in\C^{JL\times L}$ has an $L\times L$ identity matrix 
		in the first $L$ rows and zeros below, and $\vek e_{i}^{(L)}$ is the $i$th column of
		the $L\times L$ identity matrix.
		
One can derive a shifted FOM-type method by imposing the
FOM Galerkin condition on the Sylvester equation residual.
For $\widetilde{\vek D}=\vek D$, this reduces to solving a family of small shifted systems.
This was shown in \cite[Corollary 4.2]{S.1996}.  
For the $i$th shifted system, this means we 
compute $\vek t_{j}\Sf{\sigma_{i}}=\vek W_{j}\vek y_{j}\Sf{\sigma_{i}}\in\K_{j}(\vek A,\vek R_{0}^{\boldsymbol{\sigma}})$
such that 
\be\label{eqn.fom-petrov-gallerkin}
	\vek b\Sf{\sigma_{i}} - (\vek A + \sigma_{i})(\vek x_{0}\Sf{\sigma_{i}} + \vek t_{j}\Sf{\sigma_{i}\vek I})\perp \K_{j}(\vek A,\vek R_{0}^{\boldsymbol{\sigma}}),
\ee
which is equivalent to solving
\be\label{eqn.FOM-linear-system-shift}
	{\vek H}_{j}\Sf{\sigma_{i}}\vek y_{j}\Sf{\sigma_{i}} =  \vek E_{1}^{(j)}\vek S_{0}\vek e_{i}^{(L)}.
\ee

We must solve $L$ linear systems.  Each system can be solved progressively, and the 
residual norms are available at each iteration without explicit computation of the 
correction.  
One can consider this method as
a generalization of sFOM \cite{Simoncini2003a} since this is built on top of the shift invariance of the Sylvester operator
of the block Krylov subspace
which is a generalization of the shift invariance of the Krylov subspace generated by one vector.
However, we reiterate that this is not simply an extension of sFOM \cite{Simoncini2003a} to block Krylov subspaces.

The block GMRES method for Sylvester equations \cite{GJR.2002,RS.2002}.
In \cite[Equation 3.1]{GJR.2002}, it is shown that if we compute an approximate
correction to the general Sylvester equation $\CT(\vek X_{0}^{\boldsymbol\sigma} + \vek T^{\boldsymbol\sigma})=\vek B^{\boldsymbol\sigma}$ where 
$\vek T^{\boldsymbol\sigma} = \vek W_{j}\vek Y_{j}\in\K_{j}(\vek A,\vek R_{0}^{\boldsymbol\sigma})$, 
with $\vek Y_{j}\in\C^{jL\times L}$,
then the resulting residual has the form
\be\nn
	\vek R_{j}^{\boldsymbol\sigma} = \vek B^{\boldsymbol\sigma} - \CT(\vek X_{0}^{\boldsymbol\sigma} + \vek T^{\boldsymbol\sigma}_{j}) = \vek W_{j+1}\prn{\overline{\vek H}_{j}\vek Y_{j} - \vek Y_{j}\widetilde{\vek D}- \vek E_{1}^{(j)}\vek S_{0}}.
\ee
Thus $\norm{\vek R_{j}^{\boldsymbol\sigma}}_{F} = \norm{\overline{\vek H}_{j}\vek Y_{j} - \vek Y_{j}\widetilde{\vek D}- \vek E_{1}^{(j)}\vek S_{0}}_{F}$,
and solving for the GMRES minimizer of the Sylvester problem is equivalent to solving the
least-squares problem
\be\label{eqn.sylvester-min}
	\overline{\vek H}_{j}\vek Y_{j}-\vek Y_{j}\widetilde{\vek D} \approx \vek E_{1}^{(j)}\vek S_{0}.
\ee
In our case $\widetilde{\vek D} = \vek D$ is diagonal, which means that just as in the case
of \cite[Corollary 4.2]{S.1996}, the residual minimization for the full Sylvester problem
decouples into a residual minimization for each shifted system over the block Krylov subspace.
\bprop
In the case that $\widetilde{\vek D}=\vek D$ is diagonal, \eqref{eqn.sylvester-min} decouples such that
solving the least-squares problem \eqref{eqn.sylvester-min} is equivalent to solving for each column of $\vek Y_{j}$ according to
\be\label{eqn.GMRES-least-squares-shift}
	\vek y_{j}\Sf{\sigma_{i}} = \argmin{\vek y\in\C^{jL}}\norm{\vek E_{1}^{(j+1)}\vek S_{0}\vek e_{i}^{(L)} - \overline{\vek H}_{j}\Sf{\sigma_{i}}\vek y}
\ee
with $\vek Y_{j}(:,i) = \vek y_{j}\Sf{\sigma_{i}}$.
\eprop
\bproof
	One observes that \eqref{eqn.sylvester-min} can be rewritten as a tensor product equation (see, e.g., \cite{Stewart1990})
	\be\nn
		\prn{\vek I_{j}\otimes\overline{\vek H}_{j} + \vek D\otimes\vek I_{m}}{\rm vec}\prn{\vek Y_{j}}\approx {\rm vec}\prn{\vek E_{1}^{(j)}\vek S_{0}}.  
	\ee
	This leads directly to the conclusion due to the diagonal structure of $\vek D$.
\eproof

This can be restated according to
the minimum residual Petrov-Galerkin condition for the $i$th shifted system at iteration $j$, i.e.,
compute $\vek t_{j}\Sf{\sigma_{i}}=\vek W_{j}\vek y_{j}\Sf{\sigma_{i}}$ such that
\be\label{eqn.gmres-petrov-gallerkin}
	\vek b\Sf{\sigma_{i}} - (\vek A + \sigma_{i}\vek I)(\vek x_{0}\Sf{\sigma_{i}} + \vek t_{j}\Sf{\sigma_{i}})\perp (\vek A + \sigma_{i}\vek I)\K_{j}(\vek A,\vek R_{0}^{\boldsymbol{\sigma}}).
\ee
	From \eqref{eqn.blsh-arnoldi-relation} and \eqref{eqn.init-resid-relation-j}, we have that $\vek y_{j}\Sf{\sigma_{i}}$
	satisfies \eqref{eqn.GMRES-least-squares-shift}.

These least squares problems can be solved using already well-described techniques
for band upper Hessenberg matrices arising in the block Arnoldi algorithm;
see, e.g.,~\cite{GS.2005, Gutknecht2008, GS.2009}.  We must solve $L$ such least squares problems.  
We cannot simultaneously factorize them all, but
we can nonetheless efficiently solve each problem at low-cost using Householder 
reflections.
As in the block GMRES case, a progressively updated least squares residual norm
is available at each iteration, and the actual correction is only constructed at
the end of a cycle or upon convergence.  

It should be noted that this method differs 
from sGMRES \cite{Frommer1998}, because here we minimize every shifted residual.
If we begin with collinear initial residuals, we
can choose which algorithm to use; cf., Section \ref{section.collinear-handle}.  
We refer to this method as sbGMRES. 
An outline of sbGMRES and sbFOM is given in 
Algorithm \ref{alg.shbl-gmres-fom}.

In the case that the initial residuals are collinear, either of the above derived algorithms may fail or be unstable
unless procedures are in place to handle this case.  If all initial residuals are collinear, then 
one could simply choose to solve \eqref{eqn.Asigxb} using one of the methods based on the shifted invariance
\eqref{eqn.collinearity}, such as sFOM \cite{Simoncini2003a} or sGMRES \cite{Frommer1998}.
This may be preferred depending on matrix dimension and structure, which determines the cost of a block matrix-vector
product as compared to that of a single-vector matrix-vector product.
If some initial residuals are linearly dependent, one can immediately use one of the well-established 
procedures to gracefully handle this situation; see Section \ref{section.linear-dependence} and references therein.  
Here, however, we briefly explore
some other methods to transform the problem such that the initial residuals are linearly independent.  We note that
this is an exploration and that Experiment \ref{section.fourMethodCompare} demonstrates that transforming the problem
to apply one of these block methods may not be the best idea.
In Experiment
\ref{section.fourMethodCompare}, we compare the method described in Section \ref{section.gmres-handle} with
simply using sFOM \cite{Simoncini2003a} or sGMRES \cite{Frommer1998}.  Two of the techniques described in 
Section \ref{section.collinear-handle} involve generating random vectors, and in Experiments \ref{section.randBlFOMEffects}
and \ref{section.randX0Effects}, we investigate the effects of using random vectors.  

\subsection{When residuals are collinear}\label{section.collinear-handle}
In the case of collinear residuals, if we want to use either of our proposed block
methods, an initial procedure is necessary to produce non-collinear
residuals compatible with our block methods.  This can be accomplished in many different ways:
by selecting a new $\vek X_{0}$ at random so the initial residuals are unrelated,
by applying a cycle of GMRES for each shifted system over the common 
single-vector Krylov subspace they share, or by applying a cycle of block FOM in which the block
is constructed with the common residual direction as the first column and random 
vectors for the other columns.  These techniques
 produce residuals which are not collinear. We
briefly expand on the latter two ideas.  Choosing a random initial approximation requires no
explanation.  Applying a cycle of GMRES
to each shifted system is deterministic. 
The other two techniques are built on generating random vectors and 
can produce different behavior for different sets of
random vectors, with variable final outcomes; see Experiments 
\ref{section.randBlFOMEffects} and \ref{section.randX0Effects} for details.
\subsubsection{One cycle of single-vector GMRES}\label{section.gmres-handle}
In the first cycle, we can generate a single-vector
Krylov subspace (due to residual collinearity) and minimize all residuals over 
this subspace.  As long as we avoid stagnation for all right-hand sides, 
the residuals 
will not be collinear at the end of the cycle. 

\subsubsection{sbFOM with random block vectors}\label{section.blFOM-random}
We can also use a FOM iteration to obtain non-collinear residuals, but 
the iteration cannot be over the single-vector Krylov subspace.  
Shifted FOM naturally
produces collinear residuals; thus, we must do something else.
Suppose that for each shifted system, we have that the initial residual 
satisfies $\vek r_{0}\Sf{\sigma_{i}} = \beta_{i}\vek v_{1}$. 
Since all initial residuals are collinear, we can build the \emph{block} Krylov
subspace 
\be\nn
	\K_{m}(\vek A, \widetilde{\vek R}_{0}^{\boldsymbol\sigma})\mwhere\widetilde{\vek R}_{0}^{\boldsymbol\sigma} = \bbmat\vek v_{1}& \widetilde{\vek v}_{2}& \cdots & \widetilde{\vek v}_{L}\ebmat
\ee
from one normalized residual $\vek v_{1}$ (since they are all the same except for scaling) 
and some randomly generated vectors $\curl{\widetilde{\vek v}_{i}}_{i=2}^{L} $, 
similar to procedures for increasing the 
block size described in \cite{B.2000,PSS.2013,S.2014}.  
This allows us to still apply the FOM
Petrov-Galerkin condition with respect to a block Krylov subspace.  
Since the collinear residual is the first column of the block,
at the end of the cycle, for the $i$th shifted system, we now solve a problem of the form
\be\nn
	{\vek y}_{m}\Sf{\sigma_{i}} = \beta_{i}{\vek H}_{m}\Sf{\sigma_{i}}^{-1}\vek E_{1}^{(m)}\vek S_{0}\vek e_{1}^{(L)}\mwhere {\vek H}_{m}\Sf{\sigma_{i}} = \vek H_{m} + \sigma_{i}\vek I
\ee
and updating 
\be\nn
	\vek x_{m}\Sf{\sigma_{i}} = \vek x_{0}\Sf{\sigma_{i}} + \vek W_{m}\vek y_{m}\Sf{\sigma_{i}}.
\ee
After the
first cycle, the residuals are no longer collinear, and we proceed as before.

\balg
\caption{sbGMRES and sbFOM - Outline}\label{alg.shbl-gmres-fom}
\SetKwInOut{Input}{Input}\SetKwInOut{Output}{Output}
\Input{$\vek A\in \Cnn$, $\curl{\sigma_{i}}_{i=1}^{L}\subset\C$, $\vek b\Sf{\sigma_{1}},\ldots,\vek b\Sf{\sigma_{L}}\in\Cn$, $\vek x\Sf{\sigma_{1}},\ldots,\vek x\Sf{\sigma_{L}}\in\Cn$, $\varepsilon > 0$ the convergence tolerance, $m>0$ a cycle length., $m_{\rm init} > 0$ an initial cycle length}
\Output{$\vek x\Sf{\sigma_{1}},\ldots,\vek x\Sf{\sigma_{L}}\in\C^{n\times p}$ such that $\norm{\vek b\Sf{\sigma_{i}} - (\vek A + \sigma_{i}\vek I)\vek x\Sf{\sigma_{i}}}\leq \varepsilon$ for all $1\leq i\leq L$}
\For{$i=1,2,\ldots, L$}{
	$\vek r\Sf{\sigma_{i}} = \vek b\Sf{\sigma_{i}} - (\vek A + \sigma_{i}\vek I)\vek x\Sf{\sigma_{i}}$
}
\If{Initial residuals are collinear with $\vek r\Sf{\sigma_{i}}=\beta_{i}\vek r\Sf{\sigma_{1}}$}{
	Render residuals non-collinear (using a method from Section \ref{section.collinear-handle})
}
\While{$\max_{1\leq i\leq L}\curl{\norm{\vek r\Sf{\sigma_{i}}}} > \varepsilon$}{
	$\vek R^{\boldsymbol{\sigma}} \leftarrow \bbmat \vek r\Sf{\sigma_{1}} & \vek r\Sf{\sigma_{2}} & \cdots & \vek r\Sf{\sigma_{L}} \ebmat\in\C^{n\times L}$\\
	Build $\K_{m}(\vek A,\vek R^{\boldsymbol{\sigma}})$ using the block Arnoldi
	method (maintaining linear independence according to Section \ref{section.linear-dependence}), 
	generating $\vek W_{m+1}\in\C^{n\times (m+1)L}$ and 
	$\overline{\vek H}_{m}\in\C^{(m+1)L\times mL}$.\\
	\For{$i=1,2,\ldots, L$}{
		\If{GMRES}{
			$\widetilde{\vek y} \leftarrow \argmin{\vek y\in\R^{mL}}\norm{\vek E_{1}^{(m+1)}\vek S_{0}\vek e_{i} - \overline{\vek H}_{m}\Sf{\sigma_{i}}\vek y}$
		}\ElseIf{FOM}{
			$\widetilde{\vek y} \leftarrow {\vek H}_{m}\Sf{\sigma_{i}}^{-1}\vek E_{1}^{m}\vek S_{0}\vek e_{i}$
		}
		$\vek x\Sf{\sigma_{i}} \leftarrow \vek x\Sf{\sigma_{i}}  + \vek W_{m}\widetilde{\vek y}$\\
		$\vek r\Sf{\sigma_{i}} \leftarrow \vek r\Sf{\sigma_{i}}  - \vek W_{m+1}\overline{\vek H}_{m}\Sf{\sigma_{i}}\widetilde{\vek y}$
	}
	
}
\ealg
\section{Recycling for shifted systems}\label{section.shifted-recycling}
Using the Sylvester equation interpretation, we now have a method to efficiently generate a Krylov subspace and
solve shifted systems without the need to enforce a collinearity restriction on the residuals. 
Such an iteration is thus a candidate for inclusion into the subspace recycling framework.  
We focus here on developing a GMRES-type method. 
We proceed by showing that the block Krylov subspace generated by a projected Sylvester
operator has the same shift invariance as \eqref{eqn.block-invariance}.  Then we show that 
we can minimize each residual for each shifted system individually using the strategy 
described in \cite{Kilmer.deSturler.tomography.2006}.

We first observe that the invariance \eqref{eqn.block-invariance} holds even when we left-multiply by the 
orthogonal projector $\vek I - \vek P^{\CC}$.  Note, we use $\CR(\cdot)$ to denote the range of the argument.
\bprop\label{prop.proj-sylvester-invariance}
	Let $\vek V\in\C^{n\times L}$ be any block of vectors such that $\vek P^{\CC}\vek V = \vek 0$, i.e., 
	$\CR(\vek V) \perp\CC$.  Then the following invariance holds.
	\be\label{eqn.block-proj-invariance}
		\K_{j}\prn{\prn{\vek I - \vek P^{\CC}}\CT, \vek V} = \K_{j}\prn{\prn{\vek I - \vek P^{\CC}}\vek A, \vek V}.
	\ee
\eprop
\bproof
	Since $\CR(\vek V) \perp\CC$, we have the equalities
	\be\nn
		\prn{\vek I - \vek P^{\CC}}\CT(\vek V) = \prn{\vek I - \vek P^{\CC}}\prn{\vek A\vek V + \vek V\widetilde{\vek D}} = \prn{\vek I - \vek P^{\CC}}\vek A\vek V + \vek V\widetilde{\vek D}.
	\ee
	This equality is true for any block of $L$ vectors in $\CC^{\perp}$.  In particular, it holds for all blocks in
	$\CR\prn{\prn{\vek I - \vek P^{\CC}}\CT}$ since this projected operator by construction produces blocks 
	orthogonal to $\CC$.  Thus it follows that
	\be\label{eqn.proj-power-ident}
		\brac{\prn{\vek I - \vek P^{\CC}}\CT}^{j}(\vek V) = \CT_{\CC}^{j}(\vek V),
	\ee
	where we define $\CT_{\CC}:\vek X\mapsto \prn{\vek I - \vek P^{\CC}}\vek A\vek X + \vek X\widetilde{\vek D}$. We see that 
	$\CT_{\CC}$ is also a Sylvester operator which has only the trivial null space when applied to 
	$\vek V$ such that the columns of $\vek V$ are in $\CC^{\perp}$
	and it thus follows from \eqref{eqn.proj-power-ident} and
	\eqref{eqn.block-invariance} that
	\be\nn
		\K_{j}\prn{\prn{\vek I - \vek P^{\CC}}\CT, \vek V} = \K_{j}(\CT_{\CC},\vek V) = \K_{j}(\prn{\vek I - \vek P^{\CC}}\vek A, \vek V),
	\ee
	which is what we sought to prove.
\eproof
It should be noted this proposition and its proof are both generalizations of a similar proposition and proof
for one shifted system in \cite[Proposition 1]{SSX.2014}.

Thus by viewing \eqref{eqn.Asigxb} from the Sylvester equation point of view 
\eqref{eqn.shifted-sylvester}, we see from Proposition \ref{prop.proj-sylvester-invariance} that
the projected Sylvester operator is shift invariant in the sense of \eqref{eqn.block-invariance}. 
In order to minimize each shifted residual over an augmented Krylov subspace, we propose a 
hybrid approach.  
We take advantage of the invariance \eqref{eqn.block-proj-invariance} and minimize each shifted residual individually according to the strategy in \cite{Kilmer.deSturler.tomography.2006}.
 
 Properly generating $\widehat{\vek X}_{0}^{\boldsymbol\sigma}$ requires a bit of thought.
 In the case of rGMRES, one simply projects the initial residual orthogonally onto $\CC̈^{\perp}$
 and updates the initial approximation accordingly.  One can do this because the projectors
 $\vek P^{\CC}$ and $\vek P_{\vek A^{\ast}\vek A}^{\CU}$ are computable using information 
 one has on hand (namely bases for $\CU$ and $\CC$) and without solving a linear system.  
 This is not so in the case of a family of shifted systems (or in the more general
 Sylvester equation case).  This was an additional concern touched upon in \cite{SSX.2014}.
 We can certainly efficiently apply the projector $\vek I - \vek P^{\CC}$ to 
 $\vek R_{0}^{\boldsymbol\sigma}$.  However, there is no way with the information on hand
 to then update $\vek X_{0}^{\boldsymbol\sigma}$.  However, individually for each shifted system,
 one can construct an oblique projector onto $\CC^{\perp}$ such that an associated projector
 is computable without solving a linear system 
 which allows us to compute an update of the approximation just for that shifted system.
 We note that the use of oblique projections in the recycling framework for non-Hermitian systems has been
 previously advocated by Gutknecht \cite{G.2012} but in a different context than what we present below.
 We also note that these oblique projectors are applied at most once, at the beginning of a cycle.
 At each iteration, we are still able to implicitely (through a Gram-Schmidt process) apply the
 orthogonal projector $\vek I - \vek P^{\CC}$.  Thus, instabilities arising from the repeated use of oblique
 projectors does not arise here.  Also, we introduce the following subspace defiition, which is useful in the following proposition.
 \bdefin
 	Let the $k$ dimensional subspaces $\CU$ and $\CC$ have the relationship $\vek A \CU = \CC$ where for any
 	basis $\curl{\vek u_{1},\ldots,\vek u_{k}}$ of $\CU$ we denote $\curl{\vek c_{1},\ldots, \vek c_{k}}$ the basis of
 	$\CC$ such that $\vek c_{i} = \vek A\vek u_{i}$ for $1\leq i\leq k$. Then for any 
 	$\sigma\in\C\setminus\curl{0}$ we define $\CC + \sigma\CU$
 	to be the subspace spanned by $\curl{\vek c_{1}+\sigma\vek u_{1},\ldots, \vek c_{k}+\sigma\vek u_{k}}$.
 \edefin
 
\noindent  We note that this definition is well-defined and independent of bases chosen for $\CU$ and $\CC$, 
as long as they satisfy the relationship in the defintion.
 \bprop\label{prop.oblique-resid-proj}
 	Let $\CU$ and $\CC$ be defined as above, and assume we have 
 	matrices $\vek U,\vek C\in\C^{n\times k}$ such that $\CR(\vek U) = \CU$, $\vek C = \vek A\vek U$,
 	and $\vek C^{\ast}\vek C = \vek I_{k\times k}$.  For any $\sigma\in\C$, suppose 
 	we have the shifted system
 	\be\nn
 		(\vek A + \sigma\vek I)\vek x\Sf{\sigma} = \vek b\Sf{\sigma}
 	\ee
 	with initial approximation $\vek x_{0}\Sf{\sigma}$ and initial residual $\vek r_{0}\Sf{\sigma}$.
 	If $\vek Q^{\CC^{\perp}}_{\sigma}$ is the oblique projector onto $\CC + \sigma\CU$ orthogonal
 	to $\CC^{\perp}$ and $\vek Q^{\CU}_{\sigma}$ is the oblique projector onto $\CU$ orthogonal to 
 	$(\vek A + \sigma\vek I)^{\ast}\CC$, then 
 	the obliquely projected residual 
 	$\widehat{\vek r}_{0}\Sf{\sigma} = \vek r_{0}\Sf{\sigma} - \vek Q^{\CC^{\perp}}_{\sigma}\vek r_{0}\Sf{\sigma}$ has associated approximation 
 	$\widehat{\vek x}_{0}\Sf{\sigma} = {\vek x}_{0}\Sf{\sigma} + \vek Q^{\CU}_{\sigma}\prn{{\vek x}_{0}\Sf{\sigma} - \vek x\Sf{\sigma}}$,
 	and both projections are computable using only $\vek U$, $\vek C$, and $\sigma$ and with the inversion
 	of a $k\times k$ matrix.  Furthermore, we have
 	that $\widehat{\vek r}_{0}\Sf{\sigma}\perp\CC$.
 \eprop
 \bproof
 	That the resulting residual is orthogonal to $\CC$ can be seen by observing that the projector 
 	$\vek I - \vek Q^{\CC^{\perp}}_{\sigma}$ projects the residual onto $\CC^{\perp}$ orthogonal to $\CC + \sigma\CU$.
 	Using the bases we possess for $\CU$, $\CC$, and $\CC + \sigma\CU$, we can explicitly construct 
 	\be\nn
 		\vek Q^{\CC^{\perp}}_{\sigma} = \prn{\vek C + \sigma\vek U}{{\prn{\vek C^{\ast}\prn{\vek C + \sigma\vek U}}}}^{-1}\vek C^{\ast} \mand \vek Q^{\CU}_{\sigma} = \vek U{\prn{\vek C^{\ast}\prn{\vek C + \sigma\vek U}}}^{-1}\vek C^{\ast}\prn{\vek A + \sigma\vek I}
 	\ee
 	where we note that $\vek C^{\ast}\prn{\vek C + \sigma\vek U}=\vek C^{\ast}\prn{\vek A + \sigma\vek I}\CU$ and
 	that
 	\be\nn
 		\vek U{\prn{\vek C^{\ast}\prn{\vek C + \sigma\vek U}}}^{-1}\vek C^{\ast}\prn{\vek A + \sigma\vek I}\prn{\vek x\Sf{\sigma}- {\vek x}_{0}\Sf{\sigma}} = \vek U\prn{\vek C^{\ast}\prn{\vek C + \sigma\vek U}}^{-1}\vek C^{\ast}\vek r_{0}\Sf{\sigma}.
 	\ee
 	Thus these projections and updates can be computed using only matrices and vectors which are already available and
 	by inverting $\vek C^{\ast}\prn{\vek C + \sigma\vek U}$,  
 	which completes the proof.
 \eproof
We now have $\widehat{\vek X}_{0}^{\boldsymbol\sigma}$ and 
$\widehat{\vek R}_{0}^{\boldsymbol\sigma}$ such that 
$\CR\prn{\widehat{\vek R}_{0}^{\boldsymbol\sigma}}\perp\CC$.  
We can now follow \cite{Kilmer.deSturler.tomography.2006} and solve
 a small least squares problem which allows us to compute $\widehat{\vek t}_{j}\Sf{\sigma}$ and 
 $\widehat{\vek s}_{j}\Sf{\sigma}$ simultaneously, thereby minimizing the full residual of each shifted
 system over 
 $\CU + \K_{j}(\prn{\vek I - \vek P^{\CC}}\vek A,\widehat{\vek R}_{0}^{\boldsymbol\sigma})$.
 
 In the proof of Proposition \ref{prop.oblique-resid-proj}, we introduced the matrices $\vek U$ and $\vek C$, the matrices whose columns
 are bases for $\CU$ and $\CC$ such that the columns of $\vek C$ are orthonormal.  Let $\vek W_{j+1}$ and $\overline{\vek H}_{j}$ be as 
 in Section \ref{section.bl-sh-arnoldi} but now for $\K_{j}\prn{\prn{\vek I - \vek P^{\CC}}\vek A,\widehat{\vek R}^{\boldsymbol\sigma}_{0}}$.  
 Furthermore, for this projected Krylov subspace, let $\vek V_{1}$ and $\vek S_{0}$ be defined as in the unprojected case so that we
 have the QR-factorization of the initial block residual $\widehat{\vek R}_{0}^{\boldsymbol\sigma}$.
 The block version of \cite{PSS.2013} rGMRES method is built upon the augmented Arnoldi relation
 \be\nn
 	\vek A\bbmat \vek U & \vek W_{j} \ebmat = \bbmat \vek C & \vek W_{j+1} \ebmat \overline{\vek G}_{j}\mwhere \overline{\vek G}_{j} = \bbmat \vek I & \vek B_{j}\\ & \overline{\vek H}_{j} \ebmat,\mand \vek B_{j} = \vek C^{\ast}\vek A\vek W_{j},
 \ee
 and the recycled MINRES method for Hermitian problems \cite{WSG.2007} also has this property.  In \cite{Kilmer.deSturler.tomography.2006},
 the authors show that with a few modifications, a shifted Arnoldi relation can be developed which allows one to compute efficiently the
 minimum residual solution over the augmented Krylov subspace.  The main issue that must be addressed is that 
 $\prn{\vek A + \sigma\vek I}\vek U = \vek C + \sigma\vek U$ has columns forming a non-orthonormal basis of $\CC + \sigma\CU$;
 and furthermore, these columns are not orthogonal to $\K_{j}\prn{\prn{\vek I - \vek P^{\CC}}\vek A,\widehat{\vek R}^{\boldsymbol\sigma}_{0}}$.
 In \cite{Kilmer.deSturler.tomography.2006}, the QR-factorization
 \be\nn
 	\bbmat \vek W_{j+1} & \vek C & \vek U \ebmat = \bbmat \vek W_{j+1} & \vek C & \widehat{\vek U} \ebmat \bbmat \vek I & \vek 0 & \vek W_{j+1}^{\ast}\vek U \\ \vek 0 & \vek I & \vek C^{\ast}\vek U \\ \vek 0 & \vek 0 & \vek N \ebmat
 \ee
 is computed so that by defining
 \be\nn
 	\overline{\vek G}_{j}\Sf{\sigma} =\bbmat \vek I & \vek 0 & \vek W_{j+1}^{\ast}\vek U \\ \vek 0 & \vek I & \vek C^{\ast}\vek U \\ \vek 0 & \vek 0 & \vek N \ebmat \bbmat \overline{\vek H}_{j}\Sf{\sigma} & \vek 0\\ \vek B_{j} & \vek I \\ \vek 0 & \sigma \vek I  \ebmat = \bbmat \overline{\vek H}_{j}\Sf{\sigma} &  \sigma\vek W_{j+1}^{\ast}\vek U\\ \vek B_{j} & \vek I + \sigma\vek C^{\ast}\vek U\\ \vek 0 & \sigma\vek N\ebmat
 \ee
 we can derive the shifted augmented Arnoldi relation
 \be\label{shift-aug-arnoldi}
 	\prn{\vek A + \sigma\vek I}\bbmat \vek W_{j} & \vek U \ebmat = \bbmat \vek W_{j+1} & \vek C & \widehat{\vek U} \ebmat\overline{\vek G}_{j}\Sf{\sigma}.
 \ee
 This relation is used to solve the MINRES shifted residual minimization problem.  Due to the fact that the MINRES algorithm is never
 restarted, the loss of residual collinearity among the shifted residuals is not problematic.  Since we are no longer restricted by the 
 collinear residual requirement, we can use \eqref{shift-aug-arnoldi} to compute the minimum residual approximation 
 of each shifted system over the augmented block Krylov subspace efficiently in the non-Hermitian case.
\bprop
 For the $i$th shifted system from the family \eqref{eqn.Asigxb}, at iteration $j$ the minimum residual corrections $\widehat{\vek s}_{j}\Sf{\sigma_{i}}\in\CU$ 
 and $\widehat{\vek t}_{j}\Sf{\sigma_{i}}\in\K_{j}\prn{\prn{\vek I - \vek P^{\CC}}\vek A,\widehat{\vek R}^{\boldsymbol\sigma}_{0}}$ satisfying 
 \be\label{eqn.full-shift-aug-ls}
\prn{\widehat{\vek s}_{j}\Sf{\sigma_{i}},\,\widehat{\vek t}_{j}\Sf{\sigma_{i}}} = \argmin{\substack{\widehat{\vek s}\in\CU\\ \widehat{\vek t}\in\K_{j}\prn{\prn{\vek I - \vek P^{\CC}}\vek A,\widehat{\vek R}^{\boldsymbol\sigma}_{0}}}}\norm{\vek b\Sf{\sigma_{i}} - \prn{\vek A + \sigma_{i}\vek I}(\widehat{\vek x}_{0}\Sf{\sigma_{i}} + \widehat{\vek s} + \widehat{\vek t})}
 \ee
 can be computed by solving the least-squares problem
 \be\label{eqn.small-shift-aug-ls}
 	\prn{\vek y_{j}\Sf{\sigma_{i}},\,\vek z_{j}\Sf{\sigma_{i}}} = \argmin{y\in\C^{jL},\,\vek z\in\C^{k}}\norm{\vek g\Sf{\sigma_{i}} - \overline{\vek G}_{j}\Sf{\sigma_{i}}\bbmat\vek y\\\vek z\ebmat}
 \ee
 and setting $\widehat{\vek s}_{j}\Sf{\sigma_{i}} = \vek U\vek z_{j}\Sf{\sigma_{i}}$ and $\widehat{\vek t}_{j}\Sf{\sigma_{i}} = \vek W_{j}\vek y_{j}\Sf{\sigma_{i}}$
 where $\vek g\Sf{\sigma_{i}} = \widehat{\vek E}\vek S_{0}(:,i)\vek e_{i}^{(L)}$ where $\widehat{\vek E}\in\C^{(j+1)L+2k\times L}$
 has the $L\times L$ identity matrix and zeros below.
 \eprop
 \bproof
 We observe first that from the definitions of $\vek Q^{\CC^{\perp}}_{\sigma_{i}}$ and $\widehat{\vek x}_{0}\Sf{\sigma_{i}}$ that 
\be\nn
	\prn{\vek A + \sigma_{i}\vek I}\widehat{\vek x}_{0}\Sf{\sigma_{i}} = \prn{\vek A + \sigma_{i}\vek I}{\vek x}_{0}\Sf{\sigma_{i}} + \vek Q^{\CC^{\perp}}_{\sigma_{i}}\vek r_{0}\Sf{\sigma_{i}}, 
\ee 
 and for all possible pairs $\widehat{\vek s}\in\CU$ and $\widehat{\vek t}\in \K_{j}\prn{\prn{\vek I - \vek P^{\CC}}\vek A,\widehat{\vek R}^{\boldsymbol\sigma}_{0}}$, 
that $\widehat{\vek s} + \widehat{\vek t} = \bbmat \vek W_{j} & \vek U \ebmat \bbmat \vek y \\ \vek z \ebmat$ for some $\vek y\in\C^{jL}$ and $\vek z\in\C^{k}$.
Thus we can rewrite
\be\nn
\norm{\vek b\Sf{\sigma_{i}} - \prn{\vek A + \sigma_{i}\vek I}(\widehat{\vek x}_{0}\Sf{\sigma_{i}} + \widehat{\vek s} + \widehat{\vek t})} = \norm{\widehat{\vek r}_{0}\Sf{\sigma_{i}} - \prn{\vek A + \sigma_{i}\vek I}\bbmat \vek W_{j} & \vek U \ebmat \bbmat \vek y \\ \vek z \ebmat }.	
\ee
Similar to \eqref{eqn.init-resid-relation-j}, we can decompose the residual.  
Thus we have the decomposition
\be\nn
\widehat{\vek r}_{0}\Sf{\sigma_{i}} = \bbmat \vek W_{j+1} & \vek C & \widehat{\vek U} \ebmat\bbmat \vek S_{0}(:,i) \\ \vek 0_{jL} \\  \vek 0_{2k} \ebmat.
\ee
This yields the result since
\bea
	\norm{\widehat{\vek r}_{0}\Sf{\sigma_{i}} - \prn{\vek A + \sigma_{i}\vek I}\bbmat \vek W_{j} & \vek U \ebmat \bbmat \vek y \\ \vek z \ebmat } &=& \norm{\bbmat \vek W_{j+1} & \vek C & \widehat{\vek U} \ebmat\prn{\bbmat \vek S_{0}(:,i) \\ \vek 0_{jL} \\  \vek 0_{2k} \ebmat - \overline{\vek G}_{j}\Sf{\sigma_{i}}\bbmat \vek y \\ \vek z \ebmat }}\nn\\
	& = & \norm{\bbmat \vek S_{0}(:,i) \\ \vek 0_{jL} \\  \vek 0_{2k} \ebmat - \overline{\vek G}_{j}\Sf{\sigma_{i}}\bbmat \vek y \\ \vek z \ebmat },\nn
\eea
which allows us to transform \eqref{eqn.full-shift-aug-ls} to  \eqref{eqn.small-shift-aug-ls} with $\vek g\Sf{\sigma_{i}}$ so defined.
 \eproof
 We end by briefly observing that performing this minimization for each shifted system can produce residuals which are not orthogonal
 to $\CC$.  Thus, at restart, we must perform an additional projection of the residuals back into $\CC^{\perp}$.  
In a method such as rGMRES \cite{Parks.deSturler.GCRODR.2005}, this would occur naturally. 
 Here this projection incurs additional
 computational expense but can be done efficiently through a Gram-Schmidt-type process.  We present an outline of this algorithm,
 called recycled shifted block GMRES (rsbGMRES), as Algorithm \ref{alg.rbs-GMRES}.
 \balg
\caption{Recycled Block Shifted GMRES - Outline}\label{alg.rbs-GMRES}
\SetKwInOut{Input}{Input}\SetKwInOut{Output}{Output}
\Input{$\vek A\in \Cnn$, $\curl{\sigma_{i}}_{i=1}^{L}\subset\C$, $\vek b\Sf{\sigma_{1}},\ldots,\vek b\Sf{\sigma_{L}}\in\Cn$, $\vek x\Sf{\sigma_{1}},\ldots,\vek x\Sf{\sigma_{L}}\in\Cn$, $\varepsilon > 0$ the convergence tolerance, $m>0$ a cycle length, $\vek U,\vek C\in\C^{n\times k}$ with $\vek C^{\ast}\vek C=\vek I_{k\times k}$}
\Output{$\vek x\Sf{\sigma_{1}},\ldots,\vek x\Sf{\sigma_{L}}\in\C^{n\times p}$ such that $\norm{\vek b\Sf{\sigma_{i}} - (\vek A + \sigma_{i}\vek I)\vek x\Sf{\sigma_{i}}}\leq \varepsilon$ for all $1\leq i\leq L$, updated subspace $\CU$}
\For{$i=1,2,\ldots, L$}{
	$\vek r\Sf{\sigma_{i}} = \vek b\Sf{\sigma_{i}} - (\vek A + \sigma_{i}\vek I)\vek x\Sf{\sigma_{i}}$
}
\While{$\max_{1\leq i\leq L}\curl{\norm{\vek r\Sf{\sigma_{i}}}} > \varepsilon$}{
	\For{$i=1,2,\ldots, L$}{
		$\widehat{\vek r}\Sf{\sigma_{i}} = {\vek r}\Sf{\sigma_{i}} - \vek Q^{\CC^{\perp}}_{\sigma_{i}}{\vek r}\Sf{\sigma_{i}}$\\
		$\widehat{\vek x}\Sf{\sigma_{i}} = {\vek x}\Sf{\sigma_{i}} + \vek Q^{\CU}_{\sigma_{i}}\prn{\vek x_{\sigma_{i}} - {\vek x}\Sf{\sigma_{i}}}$
	}
	$\widehat{\vek R}^{\boldsymbol{\sigma}} \leftarrow \bbmat \widehat{\vek r}\Sf{\sigma_{1}} & \widehat{\vek r}\Sf{\sigma_{2}} & \cdots & \widehat{\vek r}\Sf{\sigma_{L}} \ebmat\in\C^{n\times L}$\\
	Build $\K_{m}\prn{\prn{\vek I - \vek P^{\CC}}\vek A,\widehat{\vek R}^{\boldsymbol{\sigma}}}$ using the block Arnoldi
	method, generating $\vek W_{m+1}\in\C^{n\times (m+1)L}$, 
	$\overline{\vek H}_{m}\in\C^{(m+1)L\times mL}$, and $\vek B_{j}\in\C^{k\times mL}$.\\
	\For{$i=1,2,\ldots, L$}{
		Build $\overline{\vek G}_{j}\Sf{\sigma_{i}}$ and compute $\vek g\Sf{\sigma_{i}}$\\
		Compute $\prn{\vek y_{m}\Sf{\sigma_{i}},\,\vek z_{m}\Sf{\sigma_{i}} }= \argmin{y\in\C^{jL},\,\vek z\in\C^{k}}\norm{\vek g - \overline{\vek G}_{j}\Sf{\sigma_{i}}\bbmat\vek y\\\vek z\ebmat}$\\
		Set $\widehat{\vek s}_{m}\Sf{\sigma_{i}} = \vek U\vek z_{m}\Sf{\sigma_{i}} $ and $\widehat{\vek t}_{m}\Sf{\sigma_{i}} = \vek W_{m}\vek y_{m}\Sf{\sigma_{i}}$\\
		Set $\vek x_{m} = \widehat{\vek x}\Sf{\sigma_{i}} + \widehat{\vek s}_{m}\Sf{\sigma_{i}}  + \widehat{\vek t}_{m}\Sf{\sigma_{i}} $\\
		Set $\vek r_{m} = \bbmat \vek W_{j+1} & \vek C & \widehat{\vek U} \ebmat\prn{\vek g\Sf{\sigma_{i}} - \overline{\vek G}_{j}\Sf{\sigma_{i}}\bbmat \vek y_{m}\Sf{\sigma_{i}} \\ \vek z_{m}\Sf{\sigma_{i}} \ebmat }$
	}
	Update $\CU$ (and $\CC$ if we have not converged)
}
\ealg
\section{Performance of the algorithms}\label{section.perf}
In this section, we discuss performance-related topics: stagnation and the
relationship between the sbGMRES and 
sbFOM methods, residual norms of sbGMRES 
compared to single-vector GMRES, growth of the block size due to the number of
shifts, the occurrence of linear dependence in the block Arnoldi vectors, and
selection of appropriate subspaces for recycling.
\subsection{Stagnation and the relationship of block GMRES and block FOM}
\label{section.stagnation}
The sbGMRES and sbFOM are GMRES and FOM-type methods which
are defined over the same subspace. The natural question arises: 
during a cycle, can we 
relate the approximations produced by sbGMRES and 
sbFOM  in the same way that single-vector GMRES and
FOM are related; see, e.g., \cite[Section 6.5.5]{Saad.Iter.Meth.Sparse.2003}?

In the context of the Sylvester operator interpretation, such analysis has been carried out in
\cite{RS.2002}.  One sees that as with GMRES and FOM, their Sylvester counterparts are closely related,
in that if GMRES applied to the Sylvester equations stagnates at iteration $j$, then the $j$th iteration
for FOM applied to the Sylvester equations does not exist.
One can also explore the relationships of the 
standard block GMRES and block FOM methods.
We motivate this assertion by observing that when we apply 
sbGMRES or sbFOM 
to \eqref{eqn.Asigxb}, at iteration $j$, the 
approximation $\vek  X_{j}(:,i)$ for the solution to the $i$th shifted system
is the same approximation that would be produced by applying $j$ iterations of
block GMRES or block FOM, respectively, to
\be\label{eqn.blsh.dummy-prob}
	(\vek A + \sigma_{i}\vek I)\widetilde{\vek X} = \vek R_{0}^{\boldsymbol\sigma}
\ee
with the single fixed $\sigma_{i}$ and 
initial approximation $\widetilde{\vek X}_{0}=\vek 0$ and taking\linebreak 
$\vek  X_{j}(:,i)~=~\widetilde{\vek X}_{j}(:,i)$.
Thus the behavior of sbGMRES and sbFOM can also be analyzed by considering the relationship
of block GMRES and block FOM when applied to \eqref{eqn.blsh.dummy-prob}.
Such an analysis was carried out
in a companion
work to this paper \cite{S.2014-3}, in which the occurrence of stagnation during an iteration
of block GMRES (for some or all columns of the approximation) and its relation to
block FOM is characterized.

\subsection{Comparison of block GMRES to single-vector GMRES}
The performance of block methods and comparisons to single-vector 
counterparts have been well described by many different authors; see, e.g., 
\cite{GS.2009,PSS.2013,Simoncini.Conv-Block-GMRES,S.2014,Vital1990}.  In those cases,
the analysis assumed one coefficient matrix and multiple right-hand sides.  However,
much of the convergence analysis does not specifically concern the fact that the block
Krylov subspace arises from multiple right-hand sides.  The residual 
bounds in, e.g., \cite{Simoncini.Conv-Block-GMRES}, which compare a single 
residual minimized over a block Krylov subspace to the same residual minimized
over a single-vector Krylov subspace,
are simply derived from the
fact that the block GMRES minimization is performed over a larger space.  Thus
we have the following
	\be\nn
		\norm{\vek B(:,i) - (\vek A+\sigma_{i}\vek I)\vek X_{m}(:,i)} \leq \norm{\vek B(:,i) - (\vek A+\sigma_{i}\vek I)\widehat{\vek X}_{m}(:,i)},
	\ee
	where $\vek X_{j}$ is the approximation resulting from $j$ iterations of sbGMRES and 
	$\widehat{\vek X}_{j}$ be the block approximation which
	results from applying $j$ iterations of single-vector GMRES algorithm to each shifted
	system individually.
\subsection{Block size growth with the number of shifts}
These methods allow us to solve shifted linear systems simultaneously without a 
collinearity requirement and within a subspace recycling framework 
but at a cost.  If we assume all shifted residuals are 
linearly independent, and the block Arnoldi method produces no dependent vectors,
we can then observe that the block size is dependent on
the number of shifts.  As we have stated earlier, the use of a block iteration in this
context brings with it the benefits associated to high-performance parallel computing.  
However, as the number of shifts increases, the block size also increases,
and eventually the block size will be large enough that the benefits in data-movement
efficiency will no longer outweigh the costs of the larger block size.   
Thus, we must consider what modifications can be made to accommodate this situation.
The simplest would be to choose an optimal block size $P$ and solve the
shifted systems for $P$ shifts at a time.  However, an improvement on this strategy
would be to solve $P$ systems at a time and minimize the residuals of the remaining
systems according to the strategy advocated in \cite{S.2014}.
\subsection{Linear dependence of block Arnoldi vectors}
\label{section.linear-dependence}
As with any iteration built upon a block Krylov subspace, we must address the 
possibility that during the iteration, a dependent Arnoldi vector may be produced.
As was shown in \cite{GS.2009}, the notion of the grade of a Krylov subspace extends
to the block setting.  However, unlike the single-vector case, the occurrence of a 
dependent Arnoldi vector does not indicate that the method has achieved the 
grade of the block Krylov subspace (which would imply convergence).  Many different
strategies have been suggested for gracefully handling a dependent Arnoldi vector,
see, e.g, \cite{ABFH.2000,BF.2013,D.2001,FM.1997,OLeary1980}.

We advocate replacing the dependent Arnoldi vector with a randomly 
generated vector, as in \cite{L.2003,PSS.2013,RS.2006,S.2014}.  This serves
the purpose of maintaining the block size in order to continue to realize the
data movement efficiencies associated to block methods.  
Strategies to take advantage of these block method efficiencies have been proposed in
other contexts \cite{CK.2010}.  However,
unlike \cite{S.2014}, there is no need in the nonsymmetric case to 
generate these random vectors in advance.

\subsection{Selection of recycled subspace}
\label{section.subspace-selection}
The rGMRES method \cite{Parks.deSturler.GCRODR.2005} followed from the GMRES-DR method \cite{Morgan.GMRESDR.2002}
in that the authors proposed to use harmonic Ritz vectors computed at the end of each cycle to augment at the start of 
the next cycle.  The vectors used are those associated to approximations of small eigenvalues, as it is known that smaller eigenvalues
can cause a delay in convergence, and if their influence can be damped, we can achieve accelerated convergence. 
Here we also attempt to stimulate early onset of superlinear convergence, following from the
analysis in \cite{Simoncini2005}.
Gaul and Schl\"omer
observed that in the case of MINRES, there seemed to be little difference in the amount of acceleration when Ritz vectors rather than
harmonic Ritz are used for augmentation \cite{GS.2013-arXiv}.  A criteria for the selection of an optimal subspace for recycling was 
proposed by de Sturler \cite{deSturler.GCROT.1999}.

In the case of shifted systems, one must consider that we are using one augmentation space to solve a number of shifted problems
with the same eigenvectors but different eigenvalues.   Thus computing the Ritz or harmonic Ritz vectors associated to
small approximate eigenvalues of $\vek A$ may not yield a subspace $\CU$ that damps the influence of the small eigenvalues
of $\vek A + \sigma\vek I$ for $\sigma$ large enough.  It may behoove us to rotate for which matrix we compute Ritz vectors
or compute a few vectors for multiple shifts.  Also, as observed in experiments, if we compute Ritz vectors associated
to the largest Ritz values and harmonic Ritz vectors associated to the smallest harmonic Ritz values
with respect to the unshifted matrix $\vek A$, recycling with the Ritz vectors yield far superior performance, for the set of 
matrices used in these experiments.
\section{Numerical Experiments}\label{section.num-tests}
We present experiments demonstrating the performance of
sbGMRES, sbFOM, and rsbGMRES and compare their 
performance to sFOM \cite{Simoncini2003a}, sGMRES \cite{Frommer1998}, and rGMRES.  
Unless otherwise stated Ritz vectors 
associated to the largest Ritz values with respect to the base coefficient matrix were used
for recycling in rsbGMRES.  Harmonic Ritz vectors associated to the smallest
harmonic Ritz values were used for recycling with rGMRES.
These experiments were performed on a Mac Pro with two 2.26 GHz. Quad Core Intel Xeon processors, 
12 GB of 1066 MHz DDR3 main memory running OS-X 10.10.3
using MATLAB R2014b 64-bit edition.

In all comparison experiments, we judge algorithms by iteration counts rather
than timings.  The matrices used in the experiments are relatively small.
Thus the expense of a matrix-vector product will still be dominated by FLOPS
rather than data movement. In data movement costs, it has been shown
\cite{PSS.2013} that block matrix vector products are only slightly more 
expensive than single matrix-vector products (for moderately sized blocks).
Thus for problems with millions or even billions of unknowns, the benefits
of using a block method could outweigh the costs.  In these experiments, 
we will not realize those benefits.  Furthermore, all methods were implemented in
MATLAB, and thus the overhead costs of MATLAB itself render it difficult to obtain
accurate timings for these experiments.
Thus, we compare iteration counts
in the following experiments.  To be fair, though, for Experiments
\ref{section.fourMethodCompare} and \ref{section.diffRHSSeqCompare},
for our block methods, we also show the number of iterations multiplied 
by a \emph{block matvec versus single matvec cost multiplier}.  Our chosen number of 
shifts (and therefore the block size) for these experiments is $5$.  For matrices
with the structure of those used in these experiments, multiplication by a block of 
$5$ vectors takes roughly $3.3$ times as long as when multiplying by a single vector
(when making timing measurements using compiled Trilinos libraries \cite{2011a}).  
Hence a block iteration $\times\ 3.3$ is roughly comparable to a single vector iteration,
and we use this multiplier to give a fairer comparison between the block methods discussed
in this paper and the single-vector methods against which we test them.

We performed experiments with sets of matrices coming from a 
lattice quantum chromodynamics (QCD) application.
These matrices, used in 
Experiments \ref{section.shblFOM_GMRESCompare} -- \ref{section.randX0Effects},
are from two sets of Lattice QCD matrices 
(\verb|Group 1| and \verb|Group 2|) which are of sizes $3072\times 3072$ 
and $49152\times 49152$, respectively,
available at \cite{DH.2011}.
With each such matrix $\widehat{\vek A}$ a number called $\kappa_{c}$ is provided
such that $\vek I - \widetilde{\kappa}\widehat{\vek A}$ is positive-real for all 
$0\leq \widetilde{\kappa}\leq \kappa_{c}$.  We can
equivalently state that the matrix $-\widehat{\vek A} + \kappa\vek I$ is 
positive-real for all $\frac{1}{\kappa_{c}}\leq \kappa< \infty$.
For each $\widehat{\vek A}$, we generate a base matrix
$\vek A = -\widehat{\vek A} + (10^{-3} + \kappa_{c})\vek I$ and choose
only positive shifts to create our shifted family of linear systems.  
Coefficient matrices with smaller shifts yield more poorly conditioned
systems requiring more iterations to solve.

\subsection{Convergence Curve Comparison}\label{section.shblFOM_GMRESCompare}
	
\begin{figure}[htb]
\hfill
\includegraphics[scale=0.35]{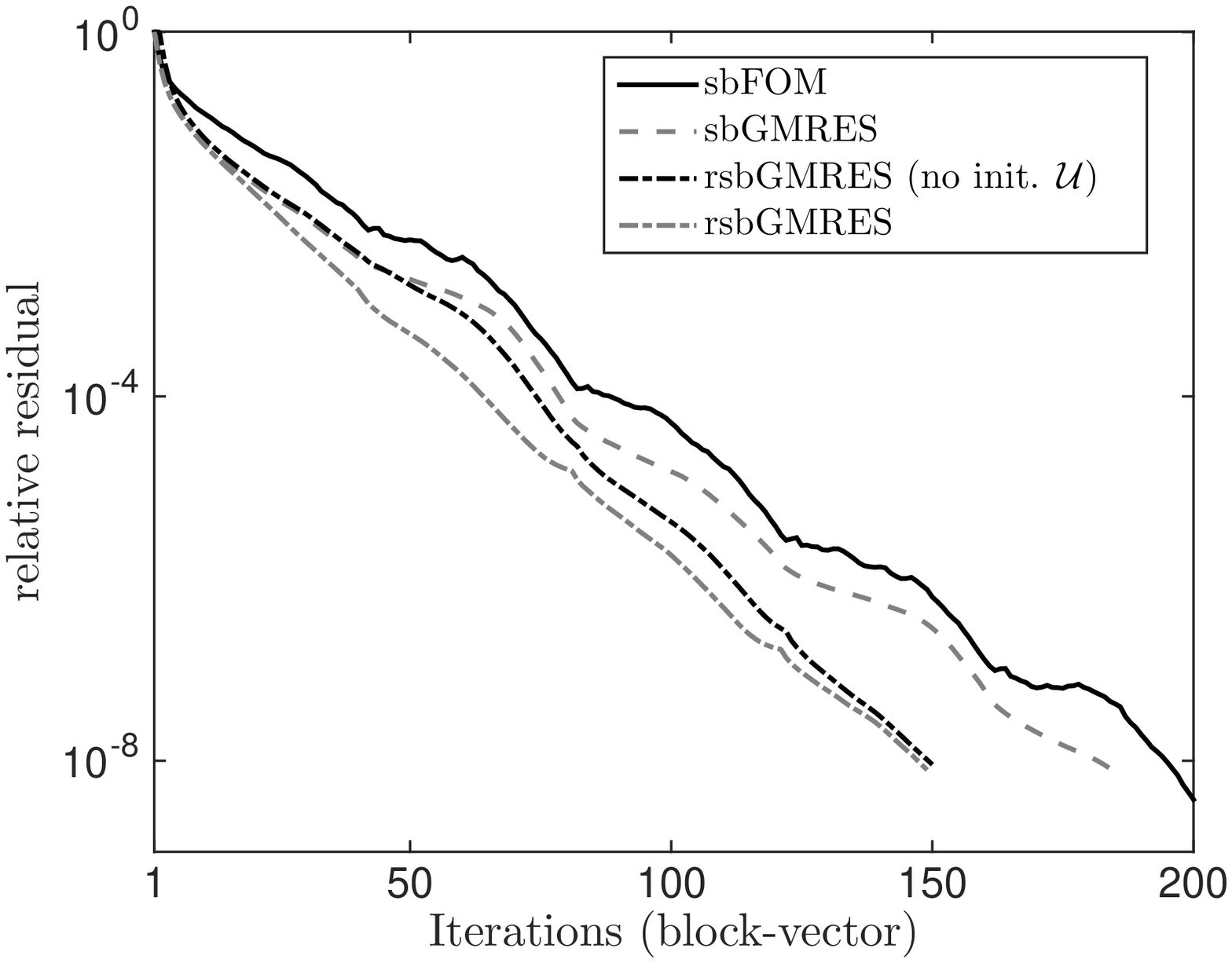}
\hfill
\begin{picture}(0,0)
\end{picture}
\caption{Convergence comparison between sbGMRES. 
sbFOM, and rsbGMRES (with and without an initial $\CU$) 
for restart parameter $m=40$, recycled subspace dimension $k=20$, and twelve shifts for the fourth matrix from {\tt Group 1}.
Initial subspace $\CU$ was generated by applying rsbGMRES to a problem with the same matrix but different right-hand side, and saving the outputted $\CU$.
We began with non-collinear residuals.
\label{figure.allMethodCompare_fig}}
\end{figure}

	In Figure \ref{figure.allMethodCompare_fig}, we compare
	the convergence histories measured with the 
	Frobenius norm of the relative residual of sbFOM, sbGMRES, and rsbGMRES. 
	The matrix used is the fourth from \verb|Group 1|, and the shifts
 	were 
 	\be\nn
 		\curl{.0001, .0002, .0003, .0004, .001, .002, .003, .004, .01, .02, .03, .04}.
 	\ee
	
 	The right-hand side for each shifted system is the same randomly
 	generated vector, and since our initial approximation for each shifted
 	system is the zero vector, we begin with noncollinear residuals. 
	Observe that initially the convergence of rsbGMRES seems to be accelerated when an
	initial recycled subspace is given, but in the end it leads to only a one iteration
	improvement.

\subsection{Mat-Vec counts for collinear initial residuals}\label{section.fourMethodCompare}

\begin{table}[htb]
\caption{Comparison of sbGMRES, sbFOM, and rsbGMRES with the sGMRES \cite{Frommer1998} and the sFOM \cite{Simoncini2003a} for four matrices from {\tt Group 2}.  For three of the matrices, the proposed methods yield improved matvec counts, but greater costs when the block matrix-vector product cost muliplier is used. This ratio will be different for different matrices and in different computing environments. The other matrices not show yielded similar results. For rsbGMRES, Ritz vectors with respect $\vek A$ associated to the largest Ritz values were recycled.}
\label{table.sixMethodCompare}
\begin{center}
\begin{tabular}{c|c|c|c|c}
Method & Matrix 1 & Matrix 2 & Matrix 3 & Matrix 4 \\
\hline
sbGMRES (mv $\times 3.3$) & 499\ (1647) & 1058\ (3491) & 442\ (1459) & 470\ (1551) \\
sGMRES \cite{Frommer1998} & 811 & 922 & 657 & 588 \\
rsbGMRES (mv $\times 3.3$) & 457\ (1508) & 443\ (1462) & 421\ (1389) & 415\ (1370) \\
sbFOM (mv $\times 3.3$) & 701\ (2313) & 869\ (2868) & 622\ (2053) & 600\ (1980) \\
sFOM \cite{Simoncini2003a} & 1138 & 1146 & 714 & 652 \\
\end{tabular}
\end{center}
\end{table}

	In Table \ref{table.sixMethodCompare}, we compared the performance
	of all five methods for 
	solving all seven systems from \verb|Group 2|.  
	In this situation, we begin with collinear right-hand sides.  For sbFOM and sbGMRES, an
	initial cycle of GMRES, as described in Section \ref{section.gmres-handle}, was 
	performed to render the initial residual noncollinear.  For rsbGMRES, this occurs naturally due to the 
	initial oblique projections of the residuals.
	Right-hand sides were generated
	as in previous experiments, and the shifts were $\curl{.0001, .0002, .01, .02}$.  In addition to comparing block iteration counts 
	versus single iteration counts, we also provide block iteration counts $\times\ 3.3$
	as described earlier.  Using this metric, the sGMRES \cite{Frommer1998} and
	sFOM \cite{Simoncini2003a} outperform sbGMRES, sbFOM, and rsbGMRES.  For all non-recycling methods, 
	cycle length $m=50$.  For rsbGMRES cycle length $m=25$, and $k=25$.

\subsection{Mat-Vec counts for unrelated initial residuals}\label{section.diffRHSSeqCompare}

\begin{table}[htb]
\caption{Comparison of block srGMRES, block sGMRES, and block sFOM with their sequentially applied counterparts when each initial residuals of the shifted systems are not collinear for one matrix from {\tt Group 2}. Both in terms of iteration counts and when multiplying block iterations by the ratio of block matvec cost to single matvec cost, we see that our proposed methods outperform their single-vector counterparts.}
\label{table.diffRHSSeqCompare_recyc}
\begin{center}
\begin{tabular}{c|c|c}
Method & Matvecs & Block Matvecs $\times 3.3$ \\
\hline
sbGMRES & 1027 & 3389 \\
rsbGMRES (Ritz) & 883 & 2914 \\
rsbGMRES (harmonic Ritz) & 1150 & 3795 \\
Sequential GMRES & 6036 & * \\
Sequential rGMRES & 5111 & * \\
sbFOM & 1252 & 4132 \\
Sequential FOM & 6647 & * \\
\end{tabular}
\end{center}
\end{table}

	In Table \ref{table.diffRHSSeqCompare_recyc}, we compare the performance of the block
	methods when the right-hand sides are unrelated, i.e., a situation in which
	sFOM \cite{Simoncini2003a} and sGMRES \cite{Frommer1998} are not applicable.  We also compare 
	both methods against simply applying GMRES, FOM, and rGMRES sequentially to 
	each shifted system.  
	The right-hand side for each system is generated randomly,
	and again the shifts were $\curl{.0001, .0002, .01, .02}$.  In terms of 
	matrix-vector product counts, the block methods are clearly superior, and
	this remains the case when we use the block iteration cost multiplier.
	For unrelated right-hand sides, it is apparent that great speedups
	can be attained.  In particular, recycling of the Ritz vectors associated to the smallest
	Ritz values yielded the greatest speedups for these problems, while harmonic Ritz recycling
	actually produced inferior results, increasing the number of iterations.  Gaul \cite{Gaul.2014-phd}
	reported in his thesis instances of inappropriate recycled subspaces slowing down the iterations.  We are 
	not sure in this case why Ritz vector recycling so greatly outperforms harmonic Ritz vector recycling.  
	For all methods rsbGMRES, 
	cycle length $m=50$.  For rsbGMRES cycle length $m=25$, and for both recycling methods $k=25$.

\subsection{Effect of random block in initial cycle of block sFOM}\label{section.randBlFOMEffects}

\begin{figure}[htb]
\hfill
\includegraphics[scale=0.35]{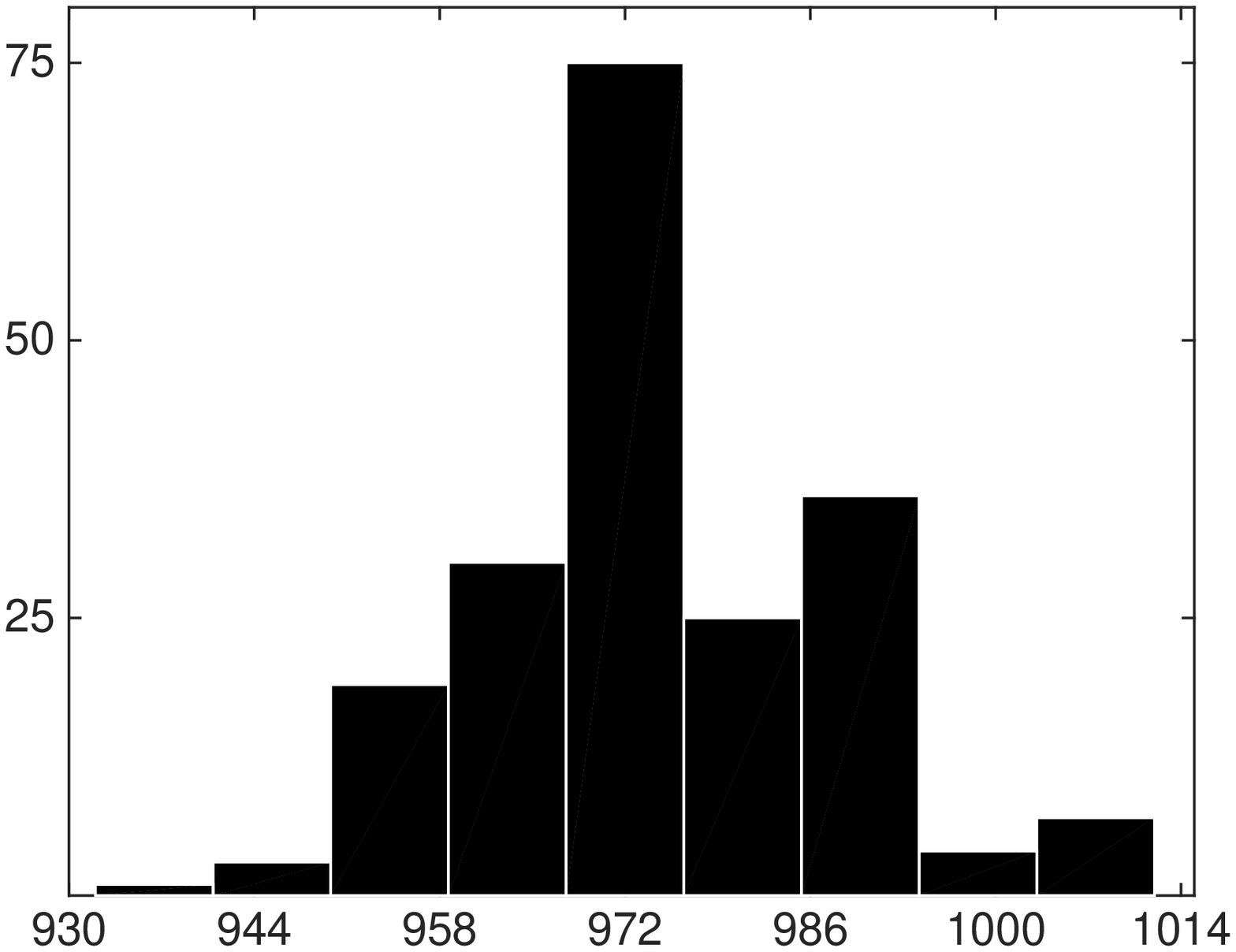}
\hfill
\begin{picture}(0,0)
\end{picture}
\caption{Histogram of matvec counts for different runs of sbFOM 
in which the initial block Krylov vector is composed of one residual and random vectors. 
Here we investigate the effect on and variation of performance for different random vectors. We have: $$mean = 974.57,\ median = 972,\ mode=972, \mand \mbox{standard deviation} = 14.344$$\label{figure.randBlFOMEffects}}
\end{figure}

	Observe that if we begin with collinear 
	residuals, we described a method in Section
	\ref{section.blFOM-random} in which we apply a block FOM cycle
	in which the block Krylov subspace is generated by one residual and $s-1$ random
	vectors.  The question arises: how variable is the performance of this method
	for different sets of random vectors.  To shed light on the answer, for
	the first matrix from \verb|Group 2|, the same shifts as in Experiment 
	\ref{section.diffRHSSeqCompare}, the same right-hand side generated
	randomly for all shifted systems, and a zero vector initial approximation
	for all systems, we applied sbFOM implemented
	with this random vector strategy to these shifted systems
	$200$ times.  Since the
	initial residuals are collinear, an initial cycle of block FOM with random
	vectors is executed in this situation. We recorded the number of iterations
	to convergence for each experiment, each with a different set of random
	vectors being generated in that initial cycle. In Figure \ref{figure.randBlFOMEffects}
	we plot a histogram for the $200$ iteration counts. As one can appreciate, 
	there is a large variation in performance ($\approx 60$ iterations), though
	a plurality of the iteration counts are clustered near the mean.

\subsection{Effect of random initial approximation $\vek X_{0}^{\boldsymbol\sigma}$}\label{section.randX0Effects}
	
\begin{figure}[htb]
\hfill
\includegraphics[scale=0.35]{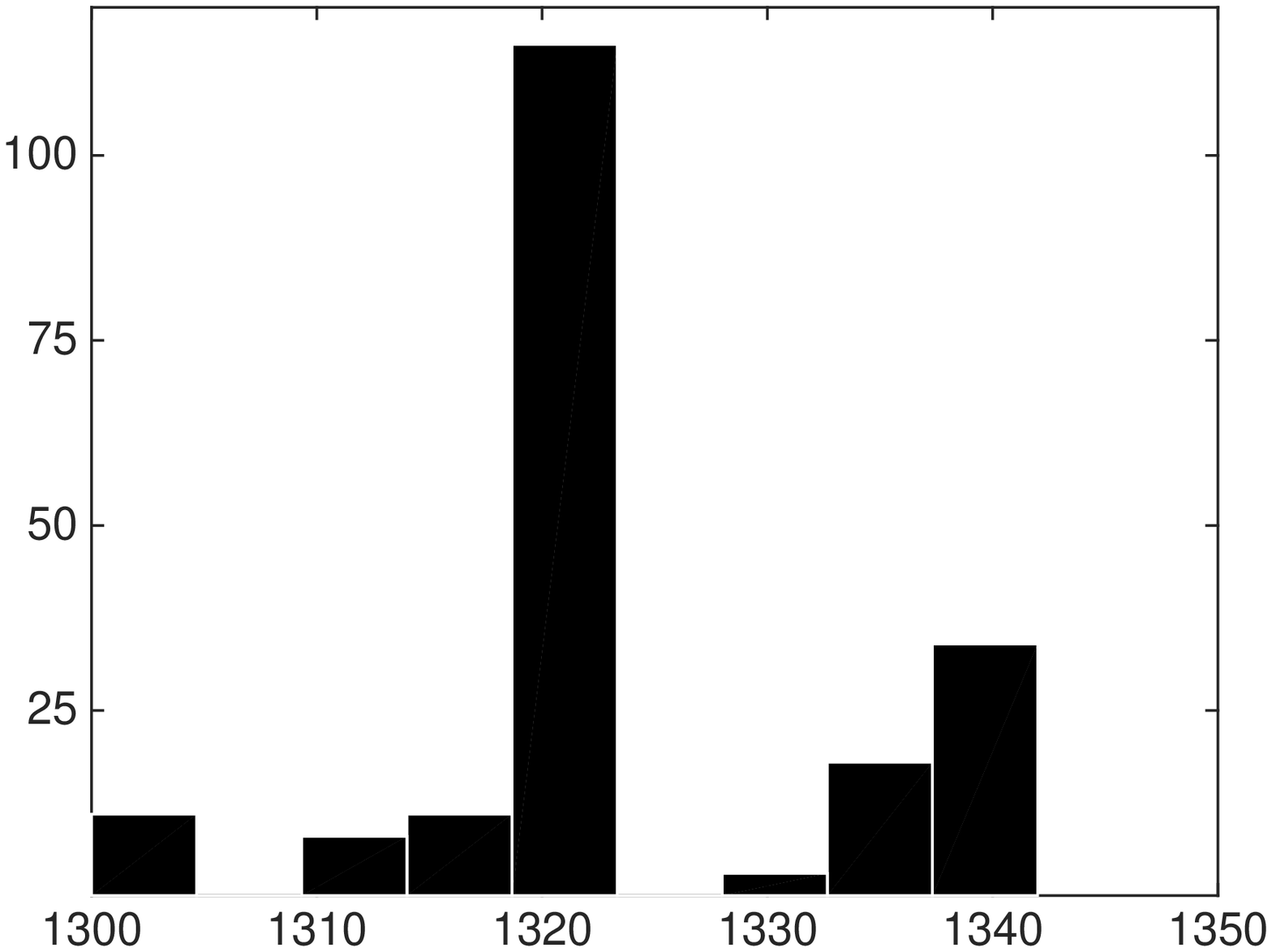}
\hfill
\begin{picture}(0,0)
\end{picture}
\caption{Histogram of matvec counts for different runs of sbFOM 
in which the initial approximation $\vek X_{0}^{\boldsymbol\sigma}$ is randomly generated. 
Here we investigate the effect on and variation of performance for different random vectors. We have: $$mean = 1324.3,\ median = 1322,\ mode=1322, \mand \mbox{standard deviation} = 9.8651$$\label{figure.randX0Effects}}
\end{figure}

	Similarly in Section \ref{section.collinear-handle}, it is mentioned that one can simply generate a random
	initial approximation in order to produce non-collinear residuals, in the case that the right-hand sides
	of \eqref{eqn.Asigxb} are collinear. Similar to Experiment \ref{section.randBlFOMEffects}, 
	for the same matrix, same shifts, and randomly generated by fixed right-hand side, we solved the shifted systems
	$200$ times using sbFOM, each time with a random starting vector, yielding a different initial residual each time.
	In Figure \ref{figure.randX0Effects}, we plot a histogram of the iteration counts produced by these $200$ experiments.
	
\subsection{Recycling when shifted system solutions depend smoothly on parameter}

\begin{figure}[htb]
\hfill
\includegraphics[scale=0.35]{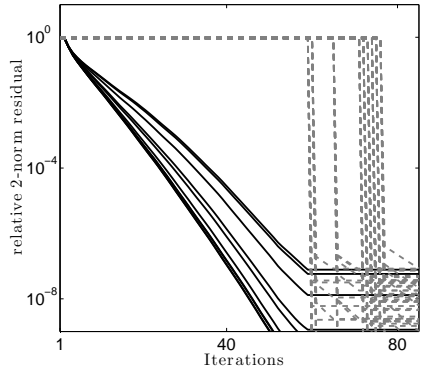}
\hfill
\begin{picture}(0,0)
\end{picture}
\caption{Here we have a family of shifted linear systems where the right-hand side (and thus the solution) depend 
smoothly on the shift.  Convergence for 100 shifts in the interval $\brac{1,2}$, where we solve only 10 systems 
and apply the theory of Kressner and Tobler by using these ten solutions as the augmentation subspace for the other 90. 
Horizontal residual curves indicates the algorithm is not acting upon those systems. 
\textbf{Black solid lines} correspond to $10$ systems first solved using sGMRES \cite{Frommer1998}, and \textbf{gray dashed lines} with 
rsbGMRES using the first $10$ solutions as augmentation vectors.  The remaining $90$ systems were solved ten-at-a-time.
\label{figure.continuousShiftDepTest}}
\end{figure}

	In this last experiment, we demonstrate that the analysis of Kressner and Tobler \cite{KT.2011} on 
	parameter dependent systems for which the dependence is sufficiently smooth offers another option for
	augmentation subspace for the rsbGMRES method.  In \cite{KT.2011}, it was shown that for a parameter dependent
	(family of) linear system(s)
	\be\nn
		\vek A\prn{\sigma}\vek x\prn{\sigma} = \vek b\prn{\sigma}
	\ee
	if the dependence on the parameter on some interval (i.e., for $\sigma\in\brac{a,b}$) is sufficiently smooth,
	then for all $\sigma$ on this interval, $\vek x\prn{\sigma}$ is well-approximated in a low-dimensional subspace.
	Thus, we can solve for a small number of parameter choices and somehow compute the rest of the solutions
	in this low-dimensional subspace.  Some specialize iterative methods for this purpose are proposed in \cite{KT.2011};
	however, one can also use this low-dimensional subspace as a high-quality augmentation subspace and use the
	rsbGMRES algorithm to solve the other systems.  In principle, this should require only initial residual projections,
	but if we compute solutions for too few parameter choices and use these solutions to augment, some further iterations
	for the remaining systems will likely be required. 
	
	We constructed an artificial example to test this idea.  We selected the first matrix from \verb|Group 1|, and 
	our interval from which we selected 100 random shifts was $\brac{1,2}$.  The shift-dependent right-hand side was constructed
	to be an infinitely smooth trigonometric function $\vek b\prn{\sigma} = \prn{\sin\prn{\sigma\frac{j\pi}{n}}}_{j}\in\C^{n}$.
	Due to the small size of $\vek A$, we can simply solve all 100 shifted systems with MATLAB's \verb|backslash| 
	and compute the dimension
	of the subspace they span.  We see that for this example, the dimension (compute numerically by taking the 
	singular value decomposition of the matrix containing solution vectors as columns) is $38$.  However, in this experiment,
	we demonstrate that one need not compute all $38$ solutions.  Here, we compute $10$ solutions for ten sample shifts, 
	using sGMRES \cite{Frommer1998}.  These solutions are used as an augmentation space, and the rest of the shifted
	systems are solved in groups of $10$.  In the results, we demonstrate that these initial ten solution vectors
	serve as an excellent augmentation space, and very few iterations are required to solve the remaining $90$
	shifted systems to tolerance.
\section{Conclusions}\label{section.conclusions}
We have shown that by taking advantage of the Sylvester equation interpretation of the
family of shifted systems, we can solve \eqref{eqn.Asigxb} with methods not restricted by residual collinearity. 
Thus we are able to propose a recycled GMRES-type method for the simultaneous solution of the systems in \eqref{eqn.Asigxb}.
This method is built upon the ``shift invariance'' of the projected Sylvester operator and the shifted system recycling
strategy of Kilmer and de Sturler \cite{Kilmer.deSturler.tomography.2006}.  
Furthermore, by basing our new methods on block Krylov subspaces,
we realize the benefits in data movement costs associated to block sparse 
matrix operations.  Numerical experiments demonstrate both the validity of the
methods and that they can outperform their single-vector counterparts.  
\section*{Acknowledgments}
We thank Jen 
Pestana for closely reading an earlier version of this manuscript and offering 
suggestions and corrections.  We also thank Martin Gutknecht and Valeria Simoncini,
who independently observed the connection between the author's initial work and the Sylvester
equation interpretation, and the author thanks Valeria Simoncini for pointers to many relevant
references in the literature.  Lastly the author thanks the referees for many helpful comments.

\bibliographystyle{siam}
\bibliography{master}
\end{document}